\documentclass{article}
\usepackage{amssymb}
\usepackage{amsmath}
\usepackage{enumerate}
\usepackage{amscd}
\usepackage{graphics}

\usepackage{hyperref}

\title{\huge Periodic Maximal surfaces in the Lorentz-Minkowski space $\l^3$}
\author{\Large Isabel Fern\'{a}ndez   \thanks{Research partially
supported by
MEC-FEDER grant number MTM2004-00160.
\newline 2000 Mathematics Subject Classification. Primary 53C50; Secondary 53C42, 53C80.
\newline Key words and phrases: maximal surfaces, periodic surfaces, conelike singularities.}
 \and \Large  Francisco J. L\'{o}pez $ ^{\ast}$  }

\usepackage{latexsym}
\usepackage{amssymb}
\usepackage{amsmath}
\usepackage{epsfig}
\usepackage[latin1]{inputenc}

\newcommand{\df}{ \stackrel{\rm def}{=}}

\def\h{\mathbb{H}}
\def\r{\mathbb{R}}
\def\n{\mathbb{N}}

\def\c{\mathbb{C}}

\def\s{\mathbb{S}}
\def\d{\mathbb{D}}
\def\l{\mathbb{L}}
\def\z{\mathbb{Z}}

\def\t{\mathbb{T}}

\def\rb{\mathcal{R}}
\def\sb{\mathcal{S}}

\def\mb{\mathcal{M}}

\def\Sg{\mathfrak{S}}
\def\pg{\mathfrak{p}}

\newenvironment{proof}{\trivlist
\item[\hskip\labelsep{\em Proof}\,:]}{\hfill{$\Box$}\endtrivlist}

\def\N{\mathcal N}
\def\D{\mathcal D}
\def\C{\mathcal C}

\def\go{G_+^\uparrow}
\def\gp{G_+}
\def\gor{G^\uparrow}
\newtheorem{lemma}{Lemma}[section]
\newtheorem{remark}{Remark}[section]
\newtheorem{theorem}{Theorem}[section]
\newtheorem{proposition}{Proposition}[section]\newtheorem{corollary}{Corollary}[section]

\newtheorem{definition}{Definition}[section]
\begin{document}
\maketitle

\begin{abstract}
A maximal surface $\sb$  with isolated singularities in a complete flat Lorentzian 3-manifold $\N$ is said to be entire if it lifts to a (periodic)  entire multigraph $\tilde{\sb}$ in $\l^3.$ In addition, $\sb$ is called of finite type if  it has finite topology, finitely many singular points and $\tilde{\sb}$ is a finitely sheeted multigraph. Complete (or  proper) maximal immersions with isolated singularities in $\N$ are entire, and entire embedded maximal surfaces in $\N$ with a finite number of singularities are of finite type.

We classify complete flat Lorentzian 3-manifolds  carrying entire maximal surfaces  of finite type, and deal with the topology, Weierstrass representation and asymptotic behavior of this kind of surfaces.   

Finally, we construct new examples of periodic entire embedded maximal surfaces in $\l^3$ with fundamental piece having finitely many singularities.
\end{abstract}
\section{Introduction}

A $3$-dimensional complete flat Lorentzian manifold $\N$ is a connected $3$-manifold with a geodesically complete flat Lorentz $(2,1)$-metric. In terms of geometric structures, it can be identified with a quotient  $\l^3/G,$ where $\l^3$ is the Lorentz Minkowski space and  $G \subset \mbox{Iso}(\l^3)$ is a discrete group of isometries acting properly and freely on $\l^3.$ Therefore, the problem of classifying complete flat Lorentzian manifolds can be  rewritten in terms of discrete groups of isometries in the Lorentz-Minkowski space. By definition, a Lorentzian manifold $\N$ is said to be orthochronous if any element of $G$ preserves the future direction.


Maximal surfaces in  $\N \equiv \l^3/G$  appear as critical points (local maxima) for the area functional associated to variations by spacelike surfaces. Among spacelike surfaces, they are characterized by the property of having zero mean curvature, and besides of their mathematical interest, they play an interesting role in classical Relativity (see \cite{marsden-tipler} for more details). From the point of view of $\l^3,$ maximal surfaces in  $\l^3/G$ correspond, up to liftings, to  periodic maximal surfaces in $\l^3$ invariant under $G.$

Calabi \cite{calabi} proved that the only complete  hypersurfaces with zero mean curvature in $\l^3$ and $\l^4$ are spacelike hyperplanes, solving the so called Bernstein-type problem in dimensions $3$ and $4.$ Cheng and Yau \cite{cheng-yau} extended this result to $\l^{n+1},$ $n \geq 4.$

A basic consequence of Calabi's theorem is the systematic arranging of  groups $G \subset \mbox{Iso}(\l^3)$ for which  $\l^3/G$ contains an  entire maximal surface $\sb.$  For topological considerations, it is natural to assume that the fundamental group of the surface is isomorphic to $G.$ Under these conditions,  $G$  contains a maximal free Abelian normal subgroup $\go$ of spacelike translations of rank one or two, and the quotient $G/\go$ is either trivial, $\z_2$ or $\z_2 \times \z_2.$ Moreover, $\sb$ must be a flat cylinder, torus, Möbius strip or Klein bottle. 

Therefore, it is meaningless to consider global
problems on maximal and {\em everywhere  regular}  surfaces in 3-dimensional complete flat Lorentzian manifolds, and for a long time, the theory of maximal surfaces dealt with the existence and regularity of solutions to the boundary value problem (see for instance Bartnik and Simon
\cite{Bar-Sim}).

\begin{figure}[htpb] \label{fig:sing}
\centering
\includegraphics[width=2cm, height=2cm]{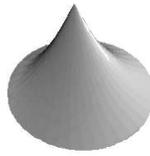} \caption{An upward pointing conelike singularity.} 	
\end{figure}

A point of a spacelike surface in $\N$ is singular if the induced metric degenerates. Although there exist spacelike surfaces admitting  curves of singularities, there are several reasons to think that {\em isolated} singular points are specially interesting. For instance,  they generally provide topological branch points of the orthogonal projection of the surface over spacelike planes. A relevant example of this kind of singularities are the so called {\em conical points} of  maximal surfaces (see Figure \ref{fig:sing}). From the PDE point of view, they appear as interior points where the elliptic operator for the maximal graph equation degenerates. Conformally, conelike singularities correspond to  analytical Jordan curves in the {\em conformal support} (see Definition \ref{def:conf}) of the maximal surface that are mapped by the immersion into  single points. Moreover, the surface reflects analytically to its  mirror about them.  From the geometrical point of view, they are points where the curvature blows up, the Gauss map has no well defined limit and the surface is locally embedded and asymptotic to the light cone. We refer to \cite{kly-mik} for a good setting.

It is then interesting  to understand the global geometry of maximal surfaces with isolated singularities.
A spacelike surface $\sb$ with singularities in $\N \equiv \l^3/G$ is said to be {\em entire} if its lifting to $\l^3$ is an entire multigraph over any spacelike plane (i.e., the orthogonal projection of $\sb$  over any spacelike plane is a branched covering). It is not hard to see that complete (or proper) spacelike surfaces in $\l^3/G$ with isolated singularities are entire, but the contrary is false. 
An entire spacelike surface $\sb$ in $\N \equiv\l^3/G$  is said to be of {\em  finite type} if it has finite topology, a finite number of singular points and it lifts to a finitely sheeted entire multigraph in $\l^3.$ The last condition means causality in the sense that timelike geodesics in $\N$ meet $\sb$ at finite time. As we will show later, entireness, properness and completenes are equivalent properties for maximal surfaces of finite type, and any entire {\em embedded} maximal surface with a finite number of singularities in $\l^3/G$ is of finite type. 

\begin{figure}[htpb]
\centering
\includegraphics[width=15cm, height=2cm]{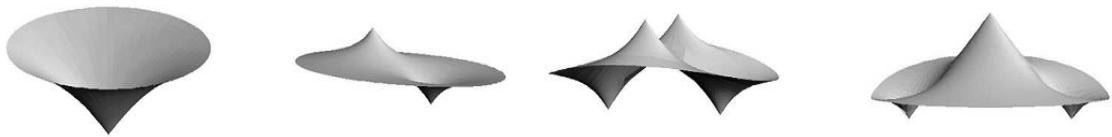}
\caption{From left to right, the catenoid, a Riemann's example and two complete maximal graphs with four and three singularities asymptotic to a plane.}\label{fig:clasicos}
\end{figure}

It is not hard to see that entire embedded maximal surfaces with isolated singularities are graphs over spacelike planes. The finite type case in $\l^3$ has been recently developed in \cite{f-l-s}. In Figure \ref{fig:clasicos} we have illustrated some known examples, including the Lorentzian catenoid and the Riemann type examples (see \cite{l-l-s} and \cite{f-l-s} for more details). 

\begin{figure}[htpb]
\centering
\includegraphics[width=5cm, height=1.5cm]{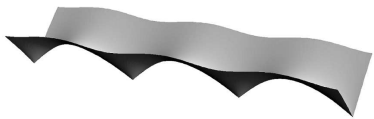}
\caption{A Scherk's type maximal surface.} \label{fig:scherk} 	
\end{figure}

\begin{figure}[htpb]
\centering
\includegraphics[width=8cm, height=2cm]{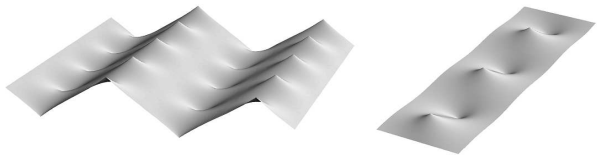}
\caption{From left to right, examples of doubly and singly periodic maximal surfaces. \label{fig:periodicas}}	
\end{figure}

However , there are many entire maximal graphs with an {\em infinite} closed discrete set of singular points prividing embedded maximal surfaces of finite type in suitable quotients $\l^3/G,$ where $G$ is a group of translations. Among them, we distinguish the one-parameter family of Scherk's type maximal surfaces (see Figure \ref{fig:scherk}), and the doubly and singly periodic examples illustrated in Figure \ref{fig:periodicas}.

The former results suggest that many open questions still remain about entire maximal surfaces with singularities.  We emphasize two of them:
\begin{enumerate}[(I)]
\item The classification of Lorentzian manifolds $\l^3/G$ containing  entire maximal surfaces of finite type, and the geometrical description of these surfaces.
\item The global geometry of entire maximal surfaces in $\l^3/G$ of finite type, starting with the periodic embedded case.
\end{enumerate}
 
The aim of this paper is to approach the above two problems.

Regarding problem (I), we will assume that the fundamental group of surfaces in $\N=\l^3/G$ surjects by the immersion onto $\Pi_1(N) \equiv G$ (replace $\N$ by a suitable covering if necessary).
If in addition this group homomorphism is injective, the surface is said to be {\em incompressible}. For instance, complete embedded maximal surfaces with a closed discrete set of singularities are always incompressible.

The classification of discrete groups $G$  acting freely and properly in $\l^3$ is far from being trivial (we refer to the surveys \cite{ch-dr-go}, \cite{zeghib} for a good setting). Among other relevant results, Margulis \cite{margulis} constructed Lorentzian manifolds $\l^3/G$ with  (non abelian) free group $G,$ and Mess \cite{mess} proved that $G$ can not be cocompact (i.e., isomorphic to the fundamental group of a closed surface with negative Euler characteristic). Mess result has played a fundamental role in the proof of the Main Theorem  below, which is devoted to classify complete flat Lorentzian manifolds $\l^3/G$ carrying entire embedded maximal surfaces of finite type.

\begin{quote}
{\bf Main Theorem: } {\em Let $\N:=\l^3/G,$ $G \neq \{\mbox{Id}\},$ be a complete flat Lorentzian manifold, and suppose that $\N$ contains an entire maximal surface $\sb$ of finite type.

Then the subgroup $\go \subset G$ of positive orthochronous isometries consists of spacelike translations and has rank $\leq 2.$ 

Moreover, $\sb$ is complete and proper, and $\sb$ is compact if and only if $\go$ has rank 2.  
If in addition $\sb$  is embedded, then it is homeomorphic to either a cylinder, a Möbius strip,  a torus or a Klein Bottle.

If $\go=G,$ the conformal support of $\sb$ is biholomorphic to a compact Riemman surface minus a finite set of interior points, and its Weierstrass data extend meromorphically to the ends of the surface.}

\end{quote}

It still remains open to know if entire  maximal surfaces in $\l^3$ with an infinite closed discrete set of singularities are complete. The Main Theorem shows that this holds when the surface is the lifting of an entire maximal surface of finite type in $\l^3/G,$ $G \neq \{\mbox{Id}\}.$

As regard to the  problem (II), we  deal with the asymptotic behavior and analytical representation of entire  maximal surfaces of finite type in $\l^3/G,$ where $G$ is a group of translations of rank $\leq 2.$ 

We obtain a geometrical formula of Jorge-Meeks type ( see \cite{jorgemeeks}, \cite{meeks-rosenberg}). These equations have significant importance, because they involve topological obstructions to the existence of this kind of surfaces. A similar formula for maxfaces in $\l^3$ was obtained by Umehara and Yamada \cite{um-ya}.

 By definition, an entire  maximal surface in $\l^3$ with a closed discrete set of singularities is said to be {\em singly periodic} (respectively, {\em doubly periodic}) if it is invariant under the action of a discrete group generated by a spacelike translation $T_1$ (respectively, two translations $T_1,$ $T_2$ spanning a group of spacelike translations), and the subsequent quotient surface is of finite type.  As a consequence of above Theorem, the lifting $\tilde{\sb} \subset \l^3$ of an entire maximal surface $\sb$ of finite type in $\N \equiv \l^3/G$ is either singly periodic or doubly periodic. Indeed, simply observe that the surface  $\tilde{\sb}/\go$ is entire and of finite type in $\l^3/\go.$

Doubly periodic entire maximal surfaces $\tilde{\sb}$ are contained in a spacelike slab invariant under $G.$ If we call $\sb_0$ the conformal support of $\sb,$ the Weierstrass representation $\Phi:=(\phi_1,\phi_2,\phi_3)$ of $\sb\equiv\tilde{\sb}/\langle T_1,T_2 \rangle$ consists of three holomorphic 1-forms  without common zeroes on the double compact Riemann surface $\Sg:=\sb_0 \cup \sb_0^*$ satisfying $J^*(\Phi_j)=-\overline{\Phi_j},$ $j=1,2,3,$ where $\sb_0^*$ is the mirror surface of $\sb_0$ and $J$ the mirror involution. In this case, $\langle T_1,T_2 \rangle$ is given by  the  real periods of $\Phi$ over ${\cal H}_1(\sb,\z).$ 

The singly periodic case admits a similar treatment. Now $\sb$ has at least two ends asymptotic (possibly with multiplicity) to flat spacelike half cylinders in $\l^3/\langle T_1 \rangle$ (exactly two Sherk's type ends in the embedded case). The Weierstrass representation $\Phi:=(\phi_1,\phi_2,\phi_3)$ of $\sb$ consists of three holomorphic 1-forms without common zeroes on $\Sg:=\sb_0 \cup \sb_0^*,$ extending meromorphically to the conformal compactification $\overline{\Sg}.$ Moreover, they have single poles at the ends of $\sb_0$ and their mirror images (ends of $\sb^*_0$), and satisfy $J^*(\Phi_j)=-\overline{\Phi_j},$ $j=1,2,3.$ Likewise, the group $\langle T_1 \rangle$ is given by the real periods of $\Phi$ over ${\cal H}_1(\sb,\z).$ The moduli problem associated to singly periodic embedded maximal surfaces has been studied in \cite{f-l-s-singly}.\\

This paper is lead out as follows:
In Section \ref{sec:prelim} we fix some notations and state some preliminary results. The proof of the Main Theorem as well as the study of the global geometry of entire maximal surfaces of finite type lies in Section \ref{sec:spacetime}. Finally, in Section \ref{sec:singly-doubly} we study the analytical representation of entire singly and doubly periodic maximal surfaces, prove the Jorge-Meeks type formula and construct the examples exhibited in Figures \ref{fig:scherk} and  \ref{fig:periodicas}.\\

{\bf Acknowledgments:} {\small We would like to thanks to M. Sanchez for helpful conversations during the preparation of this work. We are also indebted to A. Zeghib for some useful comments.} 

\section{Notations and Preliminary results} \label{sec:prelim}

We denote by $\overline{\c},$ $\d$ the extended complex plane
$\c \cup \{\infty\}$ and the unit disc $\{z \in \c \;:\; |z|<1\}.$

Throughout this paper, $\l^3$ will denote the three dimensional
Lorentz-Minkowski space $(\r^3,\langle , \rangle),$ where $\langle ,
\rangle=dx_1^2+dx_2^2-dx_3^3.$ By definition, a coordinate system $(y_1,y_2,y_3)$ in $\l^3$  is said to be a {\em $(2,1)$-coordinate system} if the Lorentzian metric is given by $dy_1^2+dy_2^2-dy_3^3.$
We say that a vector ${\bf v} \in \r^3- \{ {\bf 0} \}$ is spacelike,
timelike or lightlike if
$\|v\|^2:=\langle {\bf v}, {\bf v} \rangle$ is positive, negative or zero,
respectively. When $v$ is spacelike, $\|v\|$ is chosen non negative. The vector
${\bf 0}$ is spacelike by definition.  A plane in $\l^3$ is spacelike,
timelike or lightlike if the induced metric is Riemannian, non degenerate
 indefinite or degenerate, respectively. 

We denote by ${\cal C}_x=\{y \in \l^3 \;:\; \|y-x\|=0\}$ the lightcone with vertex at $x,$ and write $\mbox{Ext}({\cal C}_{x})=\{y \in \l^3 \;:\;\|y-x\|>0\},$ $x \in \l^3.$

We call $\h^2 = \{ (x_1,x_2,x_3) \in \r^3 \;:\;
x_1^2+x_2^2-x_3^2=-1\}$ the hyperbolic sphere in $\l^3$ of constant
intrinsic curvature $-1.$ Note that $\h^2$ has two connected  components $\h^2_+:=\h^2 \cap \{x_3 \geq 1\}$ and $\h^2_-:=\h^2 \cap \{x_3 \leq -1\}.$ The stereographic
projection $\sigma$ for $\h^2$ is defined as follows:
$$\sigma:\overline{\c} - \{|z|=1\} \longrightarrow \h^2 \,; \; z \rightarrow
\left(\frac{2 \mbox{Im} (z)}{|z|^2-1}, \frac{2 \mbox{Re} (z)}{|z|^2-1},
\frac{|z|^2+1}{|z|^2-1} \right),$$ where $\sigma(\infty)=(0,0,1).$

Denote by $\mbox{Iso}(\l^3)$ and  $\mbox{Iso}_+(\l^3)$  the group of affine isometries of $\l^3$ and  the group of positive affine isometries of $\l^3,$ respectively. By definition, an element of $\mbox{Iso}(\l^3)$ is said to be orthochronous if its associated linear  isometry  preserves $\h^2_+$ (and so $\h^2 _-$). We call $\mbox{Iso}^\uparrow(\l^3)$ and $\mbox{Iso}_+^\uparrow(\l^3)$ the group of orthochronous and orthochronous positive affine isometries, respectively.

A {\em rotation} in $\l^3$ is an element of $\mbox{Iso}_+(\l^3)$ whose fixed points are just the points on a line $\ell,$ called the axis of the rotation.  The rotation is called elliptic, hyperbolic or parabolic provided $\ell$ is a timelike, spacelike or lightlike line, respectively. 

An element of $\mbox{Iso}_+(\l^3)$ is said to be a {\em screw motion} in $\l^3$  if it is the composition of a rotation followed by a translation and is not a rotation. A screw motion is said to be elliptic, hyperbolic or parabolic provided its associated rotation is elliptic, hyperbolic or parabolic, respectively.

Positive isometries of $\l^3$ are either  translations, rotations or  screw motions. Hence, the only positive isometries in $\l^3$ without fixed points are translations of non zero vector and screw motions. Rotations and screw motions in $\l^3$  have been carefully described in the following remark.

\begin{remark}\label{re:rotations}
Given a positive isometry  $R\in\mbox{Iso}_+(\l^3),$ it is possible to find a suitable $(2,1)$-coordinate system where:
\begin{enumerate}[(1)]

\item 
$R((x,y,z))=\begin{pmatrix} \cos t& \sin t& 0\\-\sin t&\cos t&0\\0&0&1\end{pmatrix}
\begin{pmatrix} x \\ y\\z  \end{pmatrix}+  \begin{pmatrix} 0\\0\\ \lambda  \end{pmatrix},\;$
 $\lambda\in \r,\,t\in ]0,2\pi[,\;$ if $R$ is elliptic.
\item 
$R((x,y,z))=\begin{pmatrix} 1& 0&0 \\ 0& \epsilon \cosh t & \epsilon\sinh t \\ 0& \epsilon\sinh t &\epsilon\cosh t\end{pmatrix}
\begin{pmatrix} x\\y\\z \end{pmatrix}+  \begin{pmatrix} \lambda\\0\\0 \end{pmatrix},\;$
 $\lambda \in \r,\, \epsilon=\pm 1, t\in\r\;(t\neq 0 \;\rm{if}\; \epsilon=1  ),\;$
if $R$ is hyperbolic.
\item 
$R((x,y,z))=\begin{pmatrix}  1& -t&t \\ t& 1-t^2/2&  t^2/2 \\t&-t^2/2& 1+t^2/2\end{pmatrix}
\begin{pmatrix} x\\y\\z \end{pmatrix}+  \begin{pmatrix} 0\\0\\\lambda \end{pmatrix},\;$
$\lambda \in \r,\,t\neq 0,\;$ 
if $R$ is parabolic. 
\end{enumerate}
The parameter $t$ is called the {\em angle} of $R.$ Moreover, $R$ is a screw motion if and only if $\lambda \neq 0.$
\end{remark}

An elliptic or hyperbolic screw motion leaves globally (not pointwise) invariant a straight line parallel to the axis of its associated rotation, which is called the {\em axis of the screw motion}. This property does not hold for parabolic screw motions. Observe also  that any elliptic or parabolic screw motion is orthochronous, while a hyperbolic screw motion is orthochronous if and only if $\epsilon=1$ in the above representation.\\

If $\N$ is a complete flat Lorentzian 3-manifold (i.e., a 3-dimensional differential manifold endowed with a complete flat metric of index one), it is well known that the universal isometric covering of $\N$ is $\l^3$ (see for example \cite{oneill},\cite{wolf}).Thus  $\N$ can be regarded as the quotient of $\l^3$ inder the action of a discrete group $G \subset \mbox{Iso}(\l^3)$ acting freely and properly on $\l^3.$ 

If $G \subset Iso(\l^3)$ acts freely and properly on $\l^3,$  we denote by $\go=G \cap \mbox{Iso}_+^\uparrow(\l^3),$ $\gp=G \cap \mbox{Iso}_+(\l^3)$  and $\gor:=G \cap \mbox{Iso}^\uparrow(\l^3)$ the subgroup of orthochronous positive, positive  and orthochronous  isometries, respectively. We know that $\go$ (resp., $\gp,$ $\gor$) is a subgroup of index $k_+^\uparrow \leq 4$ (resp.,  $k_+ \leq 2,$ $k^\uparrow \leq 2$) in $G,$ and the  $\N_+^\uparrow:=\l^3/{\go}$ (resp., $\N_+:=\l^3/{\gp},$ $\N ^\uparrow:=\l^3/{\gor}$) is a $k_+^\uparrow$-sheeted (resp.,  $k_+$-sheeted, $k^\uparrow$-sheeted) covering of $\N:=\l^3/{G}.$  By definition, $\N_+^\uparrow$ (resp., $\N_+,$ $\N^\uparrow$) is said to be the orientable orthochronous  (resp., orientable,  orthochronous) covering of $\N.$ Note also that  $\go$ is a subgroup of index $\leq 2$ of $\gp,$ (resp., $\gor$) and $\N_+^\uparrow$ is the orthocronous (resp., orientable) covering of $\N_+$ (resp., $\N^\uparrow$).\\

Given a continuous map $X:\mb \to \N \equiv \l^3/G,$ where $\mb$ is a surface, we call $\tilde{\mb}:=\hat{\mb}/H_0,$ where $\hat{\mb}$ is the universal covering of $\mb$ and $H_0$ is the kernel of the induced group homomorphism $X_*:\Pi_1(\mb) \to \Pi_1(\N)\equiv G.$  If $\pg_1:\tilde{\mb} \to \mb$ and $\pg_2:\l^3 \to \N$ denote the natural covering projections, elementary topology implies the existence of a map $\tilde{X}:\tilde{\mb} \to \l^3$ satisfying $\pg_2 \circ\tilde{X}=X \circ \pg_1.$ This map, uniquely determined up to the initial condition, will be called a {\em lifting of $X$ to $\l^3.$} In this context, we also say that $\tilde{\mb}$ is the lifting of $\mb.$

\begin{remark} \label{re:grupo}
Thoroughout this paper, we always assume that $X_*:\Pi_1(\mb) \to \Pi_1(\N)\equiv G$ is surjective. This fact is not restrictive, because we can replace $X$ by the induced map $Y:\mb \to \l^3/X_*(\Pi_1(\mb))$ satisfying $\pg \circ Y=X,$ where $\pg:
\l^3/X_*(\Pi_1(\mb)) \to \l^3/G$ is the natural covering projection.

If in addition $X_*$ is injective,  $X$ is said to be incompressible.
\end{remark}

In what follows,  $\mb$ will denote a differential surface, $\Pi$ a spacelike plane and $\pi:\l^3 \to \Pi$ the orthogonal Lorentzian projection.

A continouos map $h:\mb \to \Pi$ is a {\em branched local homeomorphism} if given ${p} \in {\mb},$ there are open discs $V \subset{\mb}$ and $U\subset \Pi$ containing ${p}$ and $h({p}),$ respectively, such that $h(V)=U$ and $h|_{V}:V \to U$ is equivalent to the map $\varphi_n:\d \to \d,$ $\varphi_n(z)=z^n,$ $n \geq 1$ (this means that there are homeomorphisms  $\xi_1:\d \to V$ and $\xi_2:U \to \d$ such that $\xi_2 \circ (h|_{V}) \circ \xi_1=\varphi_n$). In this case, the neighbourhood $V$  is called a regular neighbourhood for $h$ at ${p}.$ The integers $n$ and $n-1$ are called the multiplicity and the branching number  of $h$ at $p,$ respectively, and $h$ is said to have a {\em branch} point at $p$ if $n-1>0.$

A branched local homeomorphism $h:\mb \to \Pi$ is said to be a {\em branched covering}  if it has the path lifting property, i.e., 
for any curve $\beta:[0,1]\to \Pi,$ $\beta(0)=x,$ and any ${p} \in h^{-1}(x),$ there is a (not necessarily unique) curve ${\alpha}:[0,1] \to\tilde{\mb}$ satisfying ${\alpha}(0)={p}$ and $h \circ {\alpha}=\beta.$ Branched coverings have well defined number of sheets.

A continous map $X:\mb\to \l^3$ is said to be an  {\em entire multigraph} over  $\Pi$ if the map $\pi \circ X:\mb \to \Pi$ is a branched covering. We also say that $X(M)$ is an entire multigraph over $\Pi.$ By definition, the number of sheets of  $X$ as multigraph over $\Pi$ is the number of sheets of $\pi \circ X.$

\subsection{Spacelike immersions with singularities}

An immersion $X:\mb \longrightarrow \N$ is said to be {\em spacelike} if for any $p \in \mb,$
the tangent plane $T_p \mb$ with the induced metric is spacelike, that is to say, the induced metric on $\mb$ is Riemannian. 
In this case, $\sb=X(\mb)$ is said to be a spacelike surface in $\N.$ 

If $\N=\l^3/G,$ where $G$ is a (possibly trivial) group of translations acting freely and properly on $\l^3,$ the locally well defined Gauss map $N_0$ of $X$ assigns to each point of $\mb$ a point of $\h^2.$  A connection argument gives that $N_0$ is globally well defined and $N_0(\mb)$ lies, up to a Lorentzian isometry, in $\h^2_-.$ This means that $\mb$ is orientable.

Let $\mb$ and  $F \subset \mb$ be  a differentiable surface and  a {\em discrete closed} subset of $\mb.$ Let $ds^2$ denote a Riemannian metric  in $\mb-F.$ Take a point $q \in F,$ an open disk ${\cal D}(q)$ in $\mb$ such that ${\cal D}(q) \cap F=\{q\}$ and  an isothermal parameter $z:{\cal D}(q)-\{q \} \to A \subset \c$ for $ds^2.$  Then write $ds^2= H |dz|^2,$ where $H(w)>0$ for any  $w \in A.$  By definition,  the Riemannian metric $ds^2$ is {\em singular} at $q$ if for any disk ${\cal D}(q)$ and any parameter $z$  as above,  the limit $\lim_{p \to q} H(z(p))$ vanishes (as a matter of fact, it suffices to check this condition just for one disc and conformal parameter). The metric $ds^2$ is said to be singular in $F$ if it is singular at any point of $F.$ In this case,  $(\mb,ds^2)$ is said to be a Riemannian surface with isolated singularities and $F$ the singular set of $(\mb,ds^2).$

\begin{definition} 
Let $\N$ be a flat  3-dimensional Lorentzian manifold, and let $X: \mb\to \N$ be a continuous map. Suppose there is a discrete and closed $F \subset \mb$ such that $X|_{\mb-F}$ is  
a spacelike immersion and $(\mb,ds^2)$ is a Riemannian surface with isolated singularities in $F,$ where $ds^2$ is the metric induced by $X.$  

Then,  $X$ is said to be  an immersion with (isolated) singularities at $F,$ and
$X(\mb)$  a spacelike surface with singularities at $X(F).$ 
\end{definition}

As stated before, throughout this paper we will always assume that the  set $F$ of  isolated singularities  is closed (and so the set of regular points is always open). This restriction is quite natural for the  global results we will approach later. Of course, this does not prevent  the existence of spacelike surfaces with points where isolated singularities accumulate, although this theory will not be treated in this paper.

In the following definition we fix the notion of {\em conformal support} of a spacelike surface with isolated singularities.
\begin{definition} \label{def:conf}
Let $X:\mb\to\N$ be an orientable spacelike immersion with a closed discrete set $F=\{q_n\;:\; n\in\Lambda\},$ $\Lambda\subset\n,$ of isolated singularities.  The surface $\mb-F$ is conformally equivalent to  $\rb-(\cup_{n \in \Lambda} D_n),$ where $\rb$ is a Riemann surface  and $\{D_n \;:\; n \in \Lambda\}$ are pairwise disjoint closed conformal discs in $\mbox{Int}(\rb)$ without accumulation in $\rb$ and whose boundaries $\{\gamma_n:=\partial(D_n)\;:\;n\in\Lambda\}$  correspond to the singularities (some discs $D_n$ could be points, and in this case $\gamma_n =D_n$).

By definition, the Riemann surface $\mb_0:=\rb\cup (\cup_{n\in\Lambda}\gamma_n)$ (with possibly non empty analytical boundary) is said to be the conformal support of the spacelike immersion $X,$ and the corresponding reparameterization $X_0:\mb_0 \to \N$ is said to be a conformal spacelike immersion.
\end{definition}
For instance, if $\mb$ is simply connected, the conformal support $\mb_0$ of $X$ is conformally equivalent to a planar circular domain (see \cite{circulardomains} for details).

A Riemann surface is said to be of {\em finite conformal type} if it is biholomorphic to a compact Riemann surface minus a finite set of interior points.

\begin{lemma}\label{lem:sing}

Let $X: \mb\to \N,$ $\Pi$ and $\pi: \l^3 \to \Pi$ be a spacelike immersion with isolated singularities,  a spacelike plane and the corresponding Lorentzian orthogonal projection, respectively.

Then, $h:=\pi\circ \tilde{X}$ is a branched local homeomorphism and its branch points correspond to the locally non embedded singularities of 
$\tilde{X},$ where $\tilde{X}:\tilde{\mb} \to \l^3$ is  the lifting of $X.$

As a consequence, if $\tilde{X}$ (or $X$) is an embedding then $\pi \circ \tilde{X}$ is a local homeomorphism.
\end{lemma}

\begin{proof}
Since $\tilde{X}$ is spacelike outside the singular set $F,$ $h$ is local homeomorphism in $\mb-F.$ Let  $p_0 \in F$ and take a compact connected neighbourhod $W \subset (\mb-F) \cup \{p_0\}$ of $p_0$ with regular boundary.  Let us see that there exists a small closed disc $U$ in $\Pi$ centered at $h(p_0)$ such that the connected component $V$ of $(h|_W)^{-1}(U)$ containing $p_0$ satisfies 
$\partial (V)\cap \partial (W)=\emptyset.$ Indeed, otherwise there exists a positive sequence $\{s_k\} \to 0$  such that, for each $k,$ $V_k \cap \partial (W)\neq \emptyset,$  where $V_k$ is the connected component of $h^{-1}(U_{k}) \cap W$ containing $p_0$ and $U_k$ is the disc of radius $s_k$ in $\Pi$ centered at $h(p_0).$ Hence,  $h^{-1}(\{p_0\})$  contains  a curve joining $p_0$ and $\partial (W),$ which contradicts that $X$ is  spacelike. 

As $\overline{V}$ is compact and $h(\partial(V))=\partial (U),$ then $h:V-h^{-1}(h(p_0)) \to U-\{h(p_0)\}$ is a proper local homeomorphism, and so a finite covering. Since $U-\{h(p_0)\}$ is a cylinder, 
we deduce that $V-h^{-1}(h(p_0))$ is a cylinder too,  $h^{-1}(h(p_0))=p_0$ and $h$ is equivalent to some $\varphi_n$ for some $n\in\n.$ 
Basic topology says that $\tilde{X}$ is an embedding locally around $p_0$ if and only if $n=1.$\end{proof}

Let $X:\mb\to\l^3$  be a spacelike immersion, and take a radial curve  $\gamma:[0,1] \to \mb$ (i.e., a curve whose image under $X$  projects orthogonally in a one to one way over a segment in some spacelike plane). If $\gamma (t)$ is non singular, $\|(X \circ \gamma)'(t)\|> 0,$ and thus integrating $\|(X \circ \gamma)(t_1)-(X \circ \gamma)(t_2)\| > 0,$ $t_1 \neq t_2.$ Therefore, if $\Pi$ is a spacelike plane and $V\subset \mb$ a regular neighborhood of $h=\pi\circ X$ such that $h(V)$ is starlike with center $h(p),$ $p \in V,$ we have: 
\begin{equation}\label{eq:cono}
X(V)-\{X(p)\}\subset\mbox{Ext}(\mathcal C_p)
\end{equation}

\begin{proposition} \label{pro:gen}
Let $X:\mb \to \l^3$  be a spacelike immersion, and suppose $h_0:=\pi_0 \circ X$ is a branched covering, where $\pi_0:\l^3 \to \Pi_0$ is the orthogonal projection over a spacelike plane $\Pi_0.$

Then, for any spacelike plane $\Pi,$ the map $h:=\pi \circ X$ is a branched covering, where $\pi:\l^3\to \Pi$ is the orthogonal projection. Furthermore, the branch points, branching numbers and number of sheets of $h$ do not depend on $\Pi.$ 
\end{proposition}
\begin{proof} Lemma \ref{lem:sing} gives that $h$ is a branched local homeomorphism. 
Pick $p\in\mb,$ take $V_0$ and $V$  regular neighbourhoods at $p$ for $h_0$ and $h,$ respectively, and  consider $\gamma \subset V_0 \cap V$  a loop winding once around $p.$ The eventual multiplicity of $p$ as branch point of $h_0$ (resp., $h$) coincides with the linking number of $X\circ \gamma$ around the timelike line $\ell_0$ (resp., $\ell$)  passing through $\tilde{X}(p)$ and orthogonal  to $\Pi_0$ (resp., $\Pi$). By Equation (\ref{eq:cono}), $p$ is a branch point of $h_0$ if and only if it is a branch point of $h,$ and with the same multiplicity.

To prove that $X$ is a multigraph over $\Pi$ it remains to check that $h:=\pi \circ X$ satisfies the path lifting property. Reasoning by contradiction, let  $\beta:[0,1]\to\Pi$ be a curve and take ${\alpha}:[0,\varepsilon[\to \mb$ a divergent inextensible lift of $\beta,$ where $\varepsilon<1.$ 

Let us see that $\lim_{t \to \varepsilon} X(\alpha(t))$ exists and is finite. Otherwise, as $\alpha$ does not extend continuously to $\varepsilon,$ we get that  either $\langle X(\alpha(t)),v \rangle$ is oscillatory or divergent, and so, $\int_{0}^{\varepsilon} |\langle (X\circ \alpha)'(t),v \rangle|dt=+\infty.$ 
Since $X$ is spacelike we infer that 
$$\infty>L(\beta)\geq \int_{0}^{\varepsilon} |\langle (X\circ \alpha)'(t),v \rangle|dt,$$
which is absurd (here $L(\cdot)$ means  {\em length} of spacelike curves with respect to the Lorentzian metric $\langle,\rangle$).
In particular, $\lim_{t \to \varepsilon} h_0(\alpha(t))$ exists, which contradicts the path lifting property for $h_0.$

Finally, assume that $h$ has a finite number $k$ of sheets and take $x \in \Pi$ such that $h^{-1}(x)$ consists of $k$ distinct points, namely, $p_1,\ldots,p_k.$ Let $V_j$ be a regular neighbourhood of $p_j$ for $h,$  and without loss of generality, assume that $h(V_j)$ is convex, $j=1,\ldots,k.$

Therefore, if $v' \in \h_+^2$ is close enough to $v,$  $X^{-1}(\ell_{v'}) \cap V_j$  consists of a unique point, $j=1,\ldots,k,$ where $\ell_{v'}=\{x+\lambda v'\;:\;\lambda\in\r\}.$ 
Let us see that that there is a neighbourhodd $W \subset \h^2_+$ of $v$ such that,  for any $v' \in W,$ $X^{-1}(\ell_{v'})$ consists of $k$ points. Otherwise, there is a sequence of vectors $\{v_n\}_{n \in \n} \subset \h^2_+$ coverging to $v$ and points 
$q_n \in X^{-1}(\ell_{v_{n}})-\cup_{j=1}^k V_j,$ $n \in \n.$ Since $\{h(q_n)\}_{n \in \n} \to x$ and $h$ is proper (recall that $h$ is  finitely sheeted), we infer that there exists a subsequence of $\{q_n \}_{n\in\n}$ converging to a point in $h^{-1}(x)=\{p_1,\ldots,p_k\},$ which leads to contradiction.

Given $W$ as above and for any $v' \in W,$  the branched covering $h':=\pi' \circ X$  has also $k$ sheets, where $\Pi'$ is a plane orthogonal to $v'$ and $\pi'$ is the orthogonal projection over $\Pi'.$
Thus,  the map applying vectors of  $\h_+^2$ into the number of sheets of the corresponding branched covering of $\mb$  is locally constant, and so globally constant ($\h_+^2$ is connected), which completes the proof.
\end{proof}

\begin{definition} \label{def:tipofinito}
A spacelike immersion $X:\mb \to \N$ with isolated singularities is said to be {\em entire} if its lifting $\tilde{X}$ is an entire multigraph over a (or any) spacelike plane. In this case, $X(\mb)$ is called an entire spacelike surface with isolated singularities.

An entire spacelike immersion $X:\mb \to \N$ with isolated singularities is said to be of finite type if
$\mb$ has finite topology, $X$  has a finite number of singularities and there is a timelike geodesic $\gamma:\r \to \N$ such that $\gamma^{-1}(X(\mb))$ is finite. We also say that $X(\mb)$ is an entire spacelike surface of finite type.
\end{definition}

The last causality condition in Definition  \ref{def:tipofinito} means that the lifting $\tilde{X}:\tilde{\mb} \to \l^3$ of $X$ is a finitely sheeted entire multigraph over spacelike planes.  
In particular, if $X:\mb \to \N$ is entire and of finite type, $\Pi$ is a spacelike plane and $\pi:\l^3 \to \Pi$ the orthogonal projection, the map $\pi \circ \tilde {X}$ is a proper map, and thus $\tilde{X}$ (and $X$) are proper.\\

A spacelike immersion $X:\mb\to \N$ with isolated singularities induces a distance function on $\mb$ in the standard way: the distance between two points is defined as the infimum length of curves joining them (note that this length is well defined even for curves passing through singular points). 
It is well known that  $\mb$  is a complete metric space if and only if every divergent path in $\mb$ has infinite length, and in this case, $X$ (resp., $X(\mb)$) is said to be a {\em complete} immersion (resp., surface).

It is natural to ask for the relationship between the concepts of properness, completeness and entireness of spacelike immersions.   Contrary to the Riemannian case, it is not hard to see that properness  does not imply completeness. 
A first approach to the proof of the following lemma in the embedded case in $\l^3$ can be found in \cite{f-l-s}.

\begin{lemma}[Global behaviour of spacelike immersions with singularities]\label{lem:graph}
Let $X:\mb\to\N$ be  spacelike immersion with isolated singularities. Then,
\begin{enumerate}[(a)]
\item  If $X$ is proper, then  $X$ is  entire.
\item If $X$ is complete, then $X$ is entire.
\item If $X$ is entire, then  $X$ is an embedding if and only if it has embedded singularities. In this case, $\tilde{X}(\tilde{\mb})$ is a graph over any spacelike plane.
\end{enumerate}
Consequently, entire embedded spacelike surfaces with finitely many  singularities are of finite type.
\end{lemma}

\begin{proof} 

To check $(a),$ first note that $X$ is proper if and only if $\tilde{X}$ is proper (use that $\Pi_1(\tilde{\mb})=\mbox{Ker}(X_*)$ and  $\N\equiv \l^3/G,$ where $G \cong \Pi_1(\mb)/\mbox{Ker}(X_*)$ acts freely and properly on $\l^3$). Hence it suffices to check that $\tilde{X}$ is an entire multigraph, provided that $\tilde{X}$ is proper. Let $\pi:\l^3 \to \Pi$ the orthogonal projection over a spacelike plane $\Pi,$ take  a smooth curve $\alpha:[0,1]\to\Pi,$  and reasoning by contradiction, suppose there is an inextensible $\tilde{\alpha}:[0,\varepsilon[\to \tilde{\mb}$ lift of $\alpha$ for $h:=\pi \circ \tilde{X},$ $\varepsilon<1.$ By the properness of $\tilde{X},$ this simply means that $\tilde{X}\circ \tilde{\alpha}$ is a divergent curve in $\r^3$. On the other hand, $\|(\tilde{X} \circ \tilde{\alpha})'(t)\|^2>0$ on non singular points implies that 
$<(\tilde{X} \circ \tilde{\alpha})',v>$ is bounded on $[0,\epsilon[,$ where $v \perp \Pi.$ This gives that
$(\tilde{X} \circ \tilde{\alpha})([0,\epsilon[)$ is bounded in $\r^3,$  which is absurd.

To prove $(b),$ take into account  that $X$ is complete if and only if $\tilde{X}$ is complete.  So, we assume that $\tilde{X}$ is complete, and reasoning again by contradiction, consider an inextensible lift $\tilde{\alpha}:[0,\varepsilon[\to \tilde{\mb}$ of $\alpha:[0,1] \to \Pi$ for $h.$ If $ds_0^2$ denotes the Riemannian metric on $\Pi$ induced by the Lorentzian one, we have $\pi^*(ds_0^2) \geq ds^2,$ where $ds^2$ is the metric in $\tilde{\mb}.$ Therefore, the length of $\tilde{\alpha}$ is less than the Riemannian length of $\alpha|_{[0,\epsilon[},$ and so it is finite. This obviously contradicts that $\tilde{\alpha}$ is divergent.

Finally, suppose all the singularities are embedded. By Lemma \ref{lem:sing}, $\pi\circ \tilde{X}$  is a local homeomorphism around embedded singularities, and thus, $h:\tilde{\mb} \to \Pi$ is a topological covering. Therefore, $h$ is a homeomorphism and $\tilde{X}$ (and $X$) an embedding, which proves $(c).$
\end{proof}

The following notion will be also useful.

\begin{definition}\label{def:pseudo}
A continuous graph $\sb\subset\l^3$ over a spacelike plane is said to be pseudo-spacelike if $\sb-\{x\}\subset\mbox{Ext}(\mathcal C_x),$ for any $x\in\sb.$
In this case, obviously $\sb$ is a graph over any spacelike plane.
\end{definition}

The next lemma provides a natural way to construct pseudo-spacelike graphs in $\l^3:$

\begin{lemma} \label{lem:sup}
Let $\sb \subset \l^3$ be an entire spacelike surface with isolated singularities, and suppose $\sb$ is a multigraph with a finite number of sheets. 
For any $v\in \h_+^2,$ call $\Pi_v:=\{z \in \l^3 \;:\; <z,v>=0\},$ and for any $z \in \Pi_v,$ denote by  $\ell_z$ the straight line $\{z+t v \;;\; t \in \r\}.$  Define $u_v^+, u_v^-:\Pi_v \to \r$ by $u_v^+(z):=\mbox{Max}\{-<p,v> \;:\; p \in  \sb \cap \ell_ z\},$  $u_v^-(z):=\mbox{Min}\{-<p,v> \;:\; p \in \sb \cap \ell_ z\}.$

Then, the  graphs $\sb_v^+:=\{z+u_v^+(z)v \;:\;  z \in \Pi_v\}$ and $\sb_v^-:=\{z+u_v^-(z)v \;:\;  z \in \Pi_v\}$ 
are pseudo-spacelike. Moreover, $\sb^+_v=\sb^+_w$ and $\sb^-_v=\sb^-_w,$ for any $v,$ $w \in \h_+^2.$

\end{lemma}
\begin{proof} 
Let $X:\mb \to \l^3$ be a  parameterizatin of $\sb,$ and as usual label $\pi_v:\l^3\to\Pi_v$ the orthogonal projection.
Let $p,q\in\sb_v^+$ and take a curve $\gamma\subset\sb_v^+$ joining $p$ and $q$ such that $\pi_v\circ\gamma$ is a segment in $\Pi_v.$ Any lift $\tilde{\gamma}$  of $\pi_v (\gamma)$ for $\pi_v \circ {X}$ is a radial curve in $\mb,$ and so from equation (\ref{eq:cono}) $X(\gamma) \subset \mbox{Ext}({\cal C}_{p}) \cup \{p\}.$ Thus $||p-q||^2>0,$ which proves the first part of the lemma.

Finally, given $x\in\sb_v^+$ and $w \in \h^2_+,$ it is easy to check that the half line $\ell$ parallel to $w,$ pointing to the future  and passing through $x$ meets $\sb_v^+$ only at $x$ (use that $\sb_v^+ -\{x\}\subset\mbox{Ext}(\mathcal C_x)$). 
This obviously implies that $x \in \sb_w^+$ and proves that $\sb_v^+=\sb_w^+.$ Likewise, $\sb_v^-=\sb_w^- .$

\end{proof}
\begin{definition} \label{def:s+}
Given an entire spacelike finitely sheeted multigraph $\sb \subset \l^3,$ we denote by $\sb^+$  (resp., $\sb^-$) the pseudo-spacelike graph $\sb_v^+ $ (resp., $\sb_v^-$), where $v $ is any timelike vector in $\h^2_+.$ 
\end{definition}

We will need some information about the asymptotic behavior of pseudo-spacelike graphs. 

\begin{lemma}\label{lem:posiciones}
Let $\sb \subset \l^3$ be an entire pseudo-spacelike graph disjoint from a lightlike line $\ell.$ Then $\sb$ is disjoint from the plane containing $\ell$ and orthogonal (in the Lorentzian sense) to $\ell.$

\end{lemma}

\begin{proof}
Without loss of generality suppose that $\ell$ is {\em above} $\sb.$ Let $\Pi$ be a spacelike plane and call $\pi:\l^3 \to \Pi$ the orthogonal projection. 
Given $q\in \ell,$ denote by $p(q):=\sb \cap \pi^{-1}(\pi(q)).$ Since $\sb-\{p(q)\} \subset \mbox{Ext}({\cal C}_{p(q)}),$ we have that $\sb \cap\mbox{Int}({\cal C}_{p(q)})^+=\emptyset,$ where $\mbox{Int}({\cal C}_{p(q)})^+$ is the future pointing convex component of $\l^3-{\cal C}_{p(q)}.$ 
Thus, $\sb$ omit the set $\cup_{q\in\ell}\mbox{Int}(\C_{p(q)})^+,$ and so 
one of the closed halfspaces determined by the plane containing $\ell$ and orthogonal in the Lorentzian sense to $\ell.$
\end{proof}


\subsection{Maximal immersions}

Let $\mb$ be a differential surface. 
A maximal immersion $X:\mb \longrightarrow \N$ is a spacelike immersion
with null mean curvature. In this case, $\sb=X(\mb)$ is said to be a maximal surface in $\N.$
Using isothermal parameters, $\mb$ can be endowed with a conformal structure. In the orientable case, $\mb$ becomes a Riemann surface. 

If $X:\mb\longrightarrow \N\equiv \l^3/G$ is a conformal maximal immersion and $G$ is a (possibly trivial) group of translations, then the Gauss map $N_0$ is well defined, $g \df \sigma^{-1} \circ N_0$ is meromorphic (here $\sigma$ is the stereographic projection), there exists a holomorphic $1$-form $\phi_3$ on $\mb$ such that 
\begin{equation}\label{eq:wei}
\phi_1=\frac{i}{2} (\frac{1}{g}-g)\phi_3, \quad \phi_2=-\frac{1}{2} (\frac{1}{g}+g) \phi_3
\end{equation}
are holomorphic on $\mb,$ $\Phi:=(\phi_1,\phi_2,\phi_3)$ never vanishes on $\mb,$ and the group generated by the real periods of $\Phi,$ i.e., $\{\mbox{Re} \big( \int_\gamma \Phi \big) \;:\; \gamma \;\mbox{closed curve in}\; \mb\},$  is contained in $G.$ Up to a translation, the immersion is given by $X= \mbox{Re}  \int_{P_0} \Phi,$ and the induced Riemannian metric on $\mb$ by
$ds^2=  |\phi_{1}|^2 +|\phi_{2}|^2- |\phi_{3}|^2 =\left( \frac{|\phi_3|}{2}
(\frac{1}{|g|}-|g|) \right)^2.$
Since $\mb$ is spacelike, then $|g| \neq 1$ on $\mb,$ and up to a Lorentzian isometry, we always assume $|g|<1.$

Conversely, if $\mb,$ $g,$ $\phi_3$ and $G$ are a Riemann surface, a holomorphic map on $\mb,$ a holomorphic $1$-form 
on $\mb$ and a group of translations acting freely and properly on $\l^3,$ respectively, satisfying that $|g| <1$  in $\mb,$ $\Phi$ is holomorphic and without zeroes in $\mb,$ and the period group generated by $\Phi$ is contained in $G,$ then $X= \mbox{Re} \big( \int
\Phi \big):\mb \to \N \equiv \l^3/G$  is a conformal maximal immersion and $N_0=\sigma \circ g$ its Gauss map. We call $(\mb, \Phi)$ (or simply
$(\mb,g,\phi_3)$)  the {\em Weierstrass representation} of $X$.
For more details see, for instance, \cite{kobayashi}.\\

Let   $X: \D\to \l^3$ be a spacelike immersed open disc with a singular point $q \in \D,$ and suppose $X|_{\D-\{q\}}$ is  maximal. There are two possibilities: either $N_0$ extend continuously to $q$ ($q$ is a {\em spacelike singular point}) or not ($q$ is a {\em lightlike} singular point). 
If $q$ is a spacelike singularity, $X$ is not a topological embedding, $\D-\{q\}$ is conformally $\d-\{0\},$ the Weierstrass data $(g,\phi_3)$ extend analytically to $q,$  $|g(q)|<1$ and $\Phi(q)=0.$ Moreover, $q$ is a branch point of $\pi \circ X$ with branching number $n_q>0,$ where $\pi$ is the orthogonal projection over a spacelike plane and $n_q$ is  the zero order of $\Phi$ at $q$ (see for instance \cite{osserman} or \cite{f-l-s}). The following lemma deals with lightlike singular points:

\begin{lemma}[Lightlike singularities \cite{kly-mik}, \cite{f-l-s}] \label{lem:sing1}
Let  $X: {\cal D}\to \l^3$ be a spacelike immersion of an open disc $\D$ with a singular point $q \in \D,$ and assume that $X|_{\D-\{q\}}$ is  a maximal immersion with a lightlike singularity at $q.$
 
Then, $\D-\{q\}$ is conformally equivalent to $\{z\in\c\; :\; 0<r<|z| < 1\},$ and
$X$ extends to a conformal map $X_0:A_r \to \l^3$ with $P_0:=X(q)=X_0(\{|z|=1\}),$ where $A_r:=\{0<r<|z|\leq 1\}.$ 
Moreover, if  $(g,\phi_3)$
are the Weierstrass  data of $X_0,$ then $|g|=1$  on $\{|z|=1\}.$
In particular, $X_0$  reflects analytically about $\{|z|=1\}$ to the mirror surface $A_r^*:=\{z\in\c\; :\; 1 \leq |z|< 1/r\},$ and so $g(J(z))=1/\overline{g}(z)$ and $J^*(\phi_3)=-\overline{\phi}_3,$ where $J(z)=1/\overline{z}$ is the mirror involution.

Furthermore, the branching number of $\pi \circ X$ at $q$ is equal to $\frac{n_q}{2}+m_q,$ where $\pi$ is the orthogonal projection over a spacelike plane, $n_q$ is the  number of zeroes counted with multiplicity of $\phi_3$ on $\{|z|=1\}$ (always even)  and $m_q$ is the degree of $g|_{\{|z|=1\}}:\{|z|=1\} \to \{|z|=1\}.$

As a consequence, $X$ is an embedding around $q$ if and only if
$g$ is injective on $\{|z|=1\}$ and $n_q=0.$ In this case and for $r_0$ close enough to $r,$ $X_0(\{0<r_0<|z| \leq 1\})$ is a graph over any spacelike plane $\Pi$ asymptotic to the top  or bottom component of ${\cal C}_{P_0},$ and $P_0$ is said to be a downward or upward pointing  conelike singularity, respectively.
\end{lemma}

Obviously, the local behavior of  singularities detailed above is valid for maximal immersions in flat Lorentzian 3-manifolds, because these manifolds are locally isometric to $\l^3.$ 
By definition, spacelike immersions with isolated singularities which are maximal outside the singular set are called {\em maximal surfaces with isolated singularities}.

>From Definition \ref{def:conf} and Lemma \ref{lem:sing1}, the boundary of the conformal support $\mb_0$ of an orientable maximal immersion $X:\mb \to \N$ with isolated singularities has as many analytical components (i.e., circles) as   lightlike singular points in $\mb.$  As in Definition \ref{def:conf},  the conformal parameterization $X_0: \mb_0 \to \N$ of $X(\mb)$ is said to be a conformal maximal immersion with isolated singularities, and the Weierstass representation $\Phi$ of $X_0$ are meromorphic data on $\mb_0$ which reflect analytically to the mirror surface $\mb_0^*.$ \\

The following definition will be crucial:

\begin{definition} \label{def:maxfinito}
A maximal immersion $X:\mb \to \N$ with isolated singularities is said to be {\em entire} if it is entire as
spacelike immersion. We also say that $X(\mb)$ is an entire maximal surface with isolated singularities.

Entire maximal surfaces of finite type are defined like in Definition \ref{def:tipofinito}. 
\end{definition}

If $X:\mb \to \N\equiv \l^3/G$ is  an entire maximal surface of finite type, any lifting $\tilde{X}:\tilde{\mb} \to \l^3$ is an entire spacelike immersion with a finite number of sheets. The converse is true, provided that $\mb\equiv\tilde{\mb}/G$ is of finite topology and $X:\mb \to \l^3/G$ is maximal outside the singular set.

Lemma \ref{lem:graph} gives that entire embedded maximal surfaces in $\l^3/G$ with a finite number of singularities are of finite type.


\section{Complete flat Lorentzian 3-manifolds and maximal surfaces} \label{sec:spacetime}

In this section we  focus  attention in the geometry of entire maximal surfaces with a finite number of singularities in complete flat Lorentzian 3-manifolds. To be more precise, we deal with maximal surfaces with a finite number of singularities and whose lifting to $\l^3$ is a finitely sheeted multigraph.
We are going to describe the non trivial subgroups $G \subset \mbox{Iso}(\l^3)$ acting freely and properly on $\l^3,$ for which $\N=\l^3/G$ admits entire maximal surfaces of finite type.
We will show  that  the index $\leq 4$ subgroup $\go$ is generated by  either one or two spacelike translations. This result allows  to control the topology and geometry of both the Lorentzian manifold and the maximal surface.
We start with the following:
 \begin{lemma}\label{lem:time+light}
There are no entire spacelike surfaces ${\sb}\subset \l^3$ of finite type invariant under  elliptic or parabolic screw motions.  
As a consequence, if ${\sb}$  is invariant under a negative isometry $R \in \mbox{Iso}_-(\l^3)$ whithout fixed points, then in a suitable $(2,1)$-coordinate system:
\begin{enumerate}[(i)]
\item Orthochronous case: $R ((x_1,x_2,x_3))=(x_1,-x_2,x_3)+(\delta,0,0), \quad \delta \neq 0,$ or
\item  Non orthochronous case: $R ((x_1,x_2,x_3))=(x_1,x_2,-x_3)+(0,\delta,0), \quad \delta \neq 0.$
\end{enumerate}
\end{lemma}

\begin{proof} Obviously,  finitely sheeted entire multigraphs are not invariant under  elliptic screw motions. Assume that $R$ is a parabolic screw motion, and up to an isometry, suppose that the axis of $\vec{R}$ is
$\ell=\{(0,s,s) \;:\; s\in \r\}$ and its translation vector is $v=(0,0,\lambda),$ $\lambda \neq 0$ (see Remark \ref{re:rotations}).

Since $R$ is orthochronous, it leaves invariant the pseudo-spacelike graph $\sb^+$ given in Definition \ref{def:s+}. Assume that ${\sb}^+\cap \ell \neq \emptyset.$ Then, for any $p\in{\sb}^+\cap \ell,$ $R(p)$ and $p$ would be points in ${\sb}$ lying on the same vertical line, which contradicts that ${\sb}^+$ is a graph over $\{x_3=0\}.$

Suppose now that ${\sb}^+\cap \ell=\emptyset.$ By Lemma \ref{lem:posiciones},  ${\sb}^+$ is contained in one the halfspaces determined by $\{(x,y,z)\in\l^3\;:\;x_3=x_2\}.$ Without loss of generality, we put ${\sb}^+\subset H=\{(x_1,x_2,x_3)\;:\;x_2 \geq x_3\}.$ 
If for any $c\in\r$ we label $H_c=\{(x_1,x_2,x_3)\in\l^3\;:\;x_2\geq x_3+c\},$ Then   ${\sb}^+=R^k({\sb}^+)\subset R^k(H)=H_{-k\lambda},$ $k\in\z,$ which is impossible.

For the final part of the lemma, note that if ${\sb}$ is invariant under a negative isometry $R,$ then $R^2$ must be an orthocronous  positive isometry different from an elliptic or parabolic screw motion. By Remark \ref{re:rotations},  $-R$ is either  a hyperbolic positive isometry, an elliptic positive isometry of angle $\pi$ or a translation. The last  possibility can not hold because $R$ has no fixed points. If $-R$ is elliptic of angle $\pi,$ $\vec{R}$ is the Lorentzian symmetry with respect to a spacelike plane. Hence $R^2$ is a {\em spacelike} translation (recall that $R^2({\sb}^+)={\sb}^+$ and ${\sb}^+$ is pseudo-spacelike), which in a suitable $(2,1)$-coordinate system leads to $(ii).$ Finally, if $-R$ is  hyperbolic, it must be non orthochronous and of angle $0$ (otherwise, $R$ would have fixed points). Reasoning as above $R^2$ is a spacelike translation, which in a suitable $(2,1)$-coordinate system corresponds to $(i).$ 
\end{proof}

The following theorem deals with the conformal type problem of entire maximal surfaces of finite type.

\begin{theorem} \label{th:parabo}
Let $\mb$ be an orientable surface, and let $X:\mb \to \N=\l^3/G,$ $G \neq \{\mbox{Id}\},$ be an entire maximal immersion  of finite type, where $G \subset \mbox{Iso}_+^\uparrow(\l^3).$ Then the conformal support $\mb_0$ of $X$ is of finite conformal type. Furthermore,
\begin{itemize}
\item If $A_0 \subset \mb$ is an annular end, the homomorphism $(X|_{A_0})_*:\Pi_1(A_0) \to G \equiv \Pi_1(\N)$ is injective,
\item If  $(X|_{A_0})_*(\Pi_1(A_0))=\langle R_0 \rangle \subset G,$ then the naturally induced immersion $Y:A_0 \to \l^3/\langle R_0 \rangle$ is an embedding outside a compact subset.
\end{itemize}
\end{theorem}
\begin{proof} First of all, we recall that $X$ and its lifting $\tilde{X}:\tilde{\mb} \to \l^3$ are proper maps (take into account that $X$ is of finite type). Denote by $\pg:\tilde{\mb} \to \mb$ the covering projection.
Since $\frac{\Pi_1(\mb)}{\mbox{Ker}(X_*)} \stackrel{X_*}{\cong} G$ (see Remark \ref{re:grupo}),   $X(\tilde{\mb})$ is invariant under $G,$ and $G$ acts freely and properly on $\tilde{\mb}$ as group of intrinsic isometries. 
>From Lemma \ref{lem:time+light} we infer that $G$ consists of  either orthochronous hyperbolic screw motions or translations, besides the identity. 

Since the conformal support $\mb_0$ of $X$ has finite topology, there is an open relatively  compact domain $D\subset\mb_0$ with analytical boundary containing $\partial (\mb_0)$ and such that $\mb_0-D$ consists of a finite family of annular ends.
Let $A_0$ be such an annular end, label $\gamma_0$ as a loop in $A_0$ generating $\Pi_1(A_0)$ and distinguish two possibilities: $\gamma_0 \in \mbox{Ker} (X_*)$ and $\gamma_0 \notin \mbox{Ker} (X_*).$\\

Let us see that the first possibility can not hold. Reasoning by contradiction, suppose $\gamma_0 \in \mbox{Ker} (X_*).$ Since $G$ is not finite, the set  $\pg^{-1}(A_0)$ is the pairwise disjoint union of infinitely many closed annular ends in $\tilde{\mb}.$ Let $A$ be any annular end in $\pg^{-1}(A_0).$  Call $\pi:\l^3\to \{x_3=0\}$ the orthogonal projection and put $h:=\pi \circ \tilde{X}.$  As
 $\tilde{X}$ is a finitely sheeted multigraph over $\{x_3=0\}$ then $h$ is a proper map. Let $U \subset \{x_3=0\}$ be an open disc containing $h(\partial(A)),$ and define $A'$ as the non compact connected component of $A \cap h^{-1} (\{x_3=0\}-U).$ Obvioulsy, $h(A') \cap (\{x_3=0\}-U)$ is open and closed in $\{x_3=0\}-U$ and so  $h(A') = \{x_3=0\}-U$. We infer that $h|_{A'}: A' \to \{x_3=0\}-U$ is a surjective proper local homeomorphism, and so a finitely sheeted covering. Likewise, given $R \in G$ and open disc $U_R \subset \{x_3=0\}$ containing $h \big(\partial(R(A))\big),$  there is an annular end $A'_R \subset R(A)$  such that $h|_{A'_R}: A'_R \to \{x_3=0\}-U_R$ is a covering. This contradicts that $h$ is finitely sheeted.\\

Thus, we can assume that $\gamma_0 \notin \mbox{Ker} (X_*).$ Let $\tilde{A}_0$  be a connected component of $\pg^{-1}(A_0) \subset \tilde{\mb},$ obviously symply connected,  and take $R_0\in G-\{\mbox{Id}\}$ such that $R_0(\tilde{A}_0)=\tilde{A}_0$ and $$\tilde{A}_0/\langle R_0\rangle\equiv A_0.$$
Up to a Lorentzian isometry, we will assume that $R_0$ leaves invariant the $x_1$-axis, that is to say, either $R_0$ is a translation of vector 
parallel to $(1,0,0)$ or ${R}_0$ is a hyperbolic screw motion with axis $\{(t,0,0) \;:\; t \in \r\}.$ This obviously implies that $\tilde{X}(\tilde{A}_0)$ can not be eventually (i.e., up to a compact subset) contained  in any halfspace orthogonal to the $x_1$-axis. 

Without loss of generality, suppose $X|_{A_0}$ and $\tilde{X}|_{\tilde{A}_0}$ are conformal parameterizations.

Since $\tilde{X}$ is spacelike, $E:=\pg (\tilde{X}^{-1}(\{x_1=c\}) \cap \tilde{A}_0)$ consists of a a family of regular analytical curves without crossing points, $c \in \r.$
Let us see that $E$ contains a divergent curve. Otherwise, any curve in $E$ must be compact, and by the maximum principle, it  bounds a compact domain in $A_0$ containing a piece of $\partial (A_0).$ Furthermore, if $X_1$ denotes the first coordinate function on $\tilde{A}_0,$ we know that $X_1 \circ R_0=X_1+\lambda,$ $\lambda \neq 0,$ and so $E$ can not contain  any homotopically non trivial loop. 
Since $\partial (A_0)$ is analytical, we infer that $E$ consists of finitely many compact curves meeting $\partial (A_0),$ and so  $\tilde{X}(\tilde{A}_0)$ is eventually  contained  in $\{x_1\geq c\}$ or $\{x_1\leq c\},$ a contradiction.

Let $\gamma \subset  \tilde{A}_0$ be a curve such that $\pg(\gamma)$ is divergent and $\tilde{X}(\gamma) \subset \{x_1=0\}.$  Up to removing a tubular neighborhood of $\partial (A_0)$ and its corresponding lifting in $\tilde{A}_0,$ we can suppose  that the initial point of $\gamma$ lies in $\partial(\tilde{A}_0).$  Let  $\Omega \subset \tilde{A}_0$ be a simply connected non compact domain given by the closure of any connected component of $\tilde{A}_0 \cap \pg^{-1}\big(A_0-\pg(\gamma)\big).$ Obviously, $\pg|_{\stackrel{\circ}{\Omega}}:{\stackrel{\circ}{\Omega}} \to A_0-\gamma$ is injective  and $\pg(\Omega)=A_0.$ Consider $\pi:\l^3\to \{x_3=0\}$  and  $h:=\pi \circ \tilde{X}$ as above.

The set $h\big(R_0^j(\partial(\Omega))\big)$ eventually bounds a half strip in $\{x_3=0\}$ parallel to the $x_2$-axis, $j \in \z.$  Let $B_j$ denote the convex hull of $h\big(\partial (R_0^j(\Omega))\big),$ $j \in \n,$ and let $\Delta_j$ denote the connected component of $\{x_3=0\}-h(\partial (R_0^j(\Omega)))$ containing  $\{x_3=0\}-\bar B_j.$ 

Let us see that $h( \Omega)$ is disjoint from $\Delta_0.$ Indeed, as $h$ is proper, then  the set $\Delta_0 \cap h( \Omega)$ is open and closed in $\Delta_0,$ and thus  either $\Delta_0 \cap h( \Omega)= \emptyset$ or $\Delta_0 \subset  h( \Omega).$ The second possibility can not hold, because otherwise we infer that $\Delta_j \subset  h( R_0^j(\Omega)),$ for any $j,$ and so  $\{x_3=0\}-\cup_{j=0}^k B_j \subset \cap_{j=0}^k h(R_0^j( \Omega)),$ for any $k >0.$ This contradicts that $h$ has a finite number of sheets.

\begin{figure}[htpb]
\centering
\includegraphics[width=13cm, height=3.5cm]{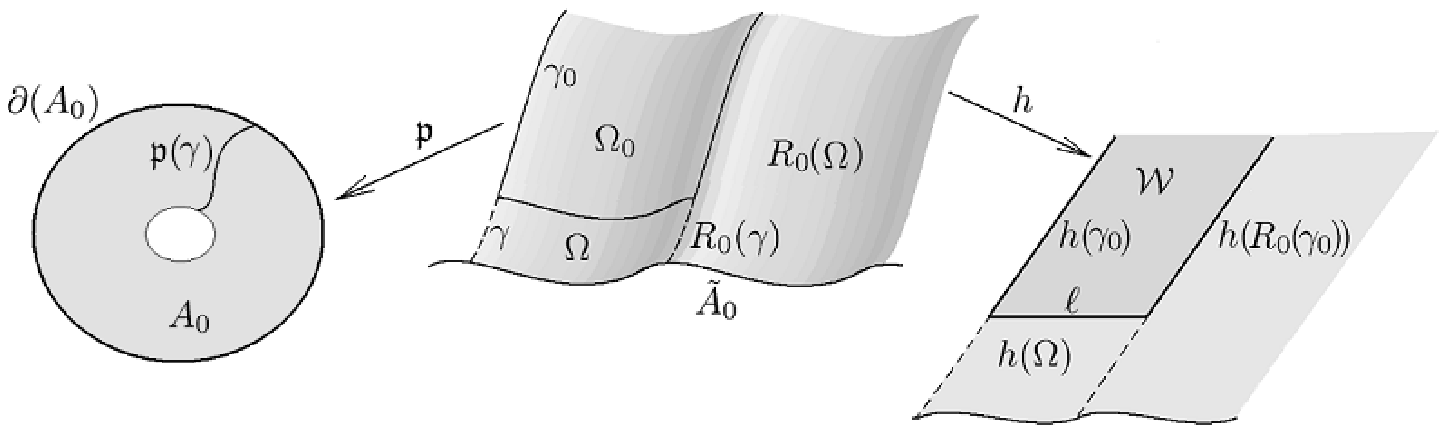}	
\end{figure}

As a consequence, $h( \Omega)$ is the union of a half strip in $\{x_3=0\}$ parallel to the $x_2$ axis and a compact set of this plane.
In the sequel, and without loss of generality, we suppose that $\Omega$ contains $\gamma \cup R_0(\gamma)$ in its boundary.  Consider a closed half strip ${\cal W} \subset h(\Omega)$ whose boundary consists of $h(\gamma_0) \cup h(R_0(\gamma_0)) \cup \ell,$ where $\gamma_0 \subset \gamma$ is a divergent curve and $\ell$ is a segment joining the initial points of 
the half lines $h(\gamma_0)$ and $h(R_0(\gamma_0))$ and contained (up to its initial and final points) in $h(\Omega)-\partial(h(\Omega)).$ Since the closure of $h(\Omega)-{\cal W}$ is compact and $h$ is proper, the maximum principle gives that the closure of $\Omega-\Omega_0$ is a compact simply connected domain with boundary $\big(\partial(\Omega) \cup h^{-1}(\ell) \big)-(\gamma_0 \cup R_0(\gamma_0)),$ where $\Omega_0:=(h|_\Omega)^{-1}({\cal W}).$  Furthermore, $h|_{\Omega_0}:\Omega_0 \to {\cal W}$ is a proper local homeomorphism, and so a covering. Since ${\cal W}$ is simply connected, we infer that $h|_{\Omega_0}$ is one to one, that is to say, $\tilde{X}(\Omega_0)$ is a graph over ${\cal W}.$

Since   $A_0-\pg(\Omega_0)$ is a tubular neighborhood of $\partial (A_0)$ with compact closure,  we can relabel $A_0 \equiv \pg(\Omega_0)$ and $\Omega \equiv \Omega_0.$ In particular, we can assume that $\tilde{X}(\Omega)$ is a graph over ${\cal W}$  and the induced map $Y:A_0 \to \l^3/\langle R_0 \rangle$ is a proper embedding.

On the other hand, it is obvious that the second coordinate function $X_2:=x_2 \circ \tilde{X}$ is bounded above or below, proper and harmonic on $\Omega.$ This gives that $\Omega$ is relative parabolic with the conformal structure induced by $\tilde{X}$ (see for instance  \cite{ahlfors}, \cite{grigor}). Let $F:\tilde{A}_0 \to \c$ be the holomorphic map given by $F:=\int \phi_1,$ where $\phi_1$ is the first 1-form in the Weierstrass representation of $\tilde{X}|_{\tilde{A}_0}.$
Observe that $F(p)=X_1(p)+i X_1^*(p),$ where $X_1=x_1 \circ \tilde{X}$ and $X_1^*$ is the harmonic conjugate of $X_1.$ 
Since $\tilde{X}(\Omega)$ is a graph over ${\cal W},$ $(h|_\Omega)^{-1}(\{x_1=c\})$ is a regular divergent arc, for any $c \in I \subset \r,$ where $I$ is a suitable compact interval of length the width of ${\cal W}.$ As a consequence, $X_1^*$ is injective on the level curves of $X_1,$ and so $F|_\Omega$ is injective.  We infer that $F(\Omega)$ is a relative parabolic domain in $\c,$ which  obviously implies that $X_1^*|_\Omega:\Omega \to \r$ is a proper harmonic map bounded above or below. 
Since $R_0^*(\phi_1)=\phi_1,$ then $F \circ R_0=F+\omega,$ where $\omega \in \c-\{0\},$ and so the map $H:=e^{\frac{2 \pi i F}{\omega}}:{A}_0 \to \c$ is well defined.  Furthermore, $\log{|H|}:A_0 \to \r$ is a harmonic proper  function bounded above or below, which proves that $A_0$ is biholomorphic to $\bar \d-\{0\}$ and concludes the proof.
\end{proof}

\begin{lemma}\label{lem:traslaciones}
Let ${\sb}\subset \l^3/G$ be an entire spacelike surface of finite type, where $G \subset \mbox{Iso}_+^\uparrow(\l^3)$ and contains a translation $T \neq \mbox{Id}.$ 
Then $G$ is a group of spacelike translations of rank at most two.
\end{lemma}

\begin{proof}
Label $\tilde{\sb}$ as a lifting of $\sb$ and write $T(x)=x+w,$ $w \neq 0.$ Since $\tilde{\sb}^+$ is a pseudo-spacelike graph invariant under $T,$  $w=T(p)-p$ is spacelike, $p \in \tilde{\sb}^+.$   

Reasoning by contradiction, suppose  $G$ contains a Lorentzian screw motion $R.$ By Lemma \ref{lem:time+light},  $R(x)=\vec{R}(x)+v,$ where $\vec{R}$ is a linear orthochronous hyperbolic rotation of non zero angle, and up to a Lorentzian isometry, $v$ lies in the axis $\ell$ of $\vec{R}.$

Define $T'=R^{-1} \circ T \circ R \in G$ and observe that $T'(x)=x+w',$ where $w'=\vec{R}^{-1}(w).$ 

Let us see that the vectors $w$ and $w'$ are linearly independent. Otherwise, $\vec{R}(w)=\lambda w,$ $\lambda \in \r.$  We infer that $w$ is a spacelike eigenvector of $\vec{R},$ that is to say, $w$ lies in the axis of $\vec{R}.$ Therefore $\lambda=1,$ $T=T'$ and $w=\mu v,$ $\mu \in \r.$  Since  $G'=\langle T,R\rangle$ acts freely and properly as a group of translations on $\ell,$ it must be cyclic. Hence there are $n,m \in \z-\{0\}$ such that $T^m \circ R^n$ fixes the origin, which contradicts that $G$ acts freely and properly on $\l^3.$

Denote by $\Pi$ the spacelike plane generated by  $w$ and $w',$ and call $G_1=\langle T,T'\rangle.$ As $\tilde{\sb}^+$ is a graph over $\Pi,$ then $\tilde{\sb}^+/G_1$ is a topological torus 
and the same holds for $\tilde{\sb}^+/G$ since it is  covered by $\tilde{\sb}^+/G_1.$ This proves that $G$ is a rank two Abelian group.

Since $G_1$ is a finite index subgroup of $G,$ we infer that $R^n \in G_1,$ for a suitable $n>0,$ and so $R^n$ is a translation. This is obviously absurd, and proves that $G$ is a translational group.

If $G$ contains two independent translations, we repeat the last argument to get that $\tilde{\sb}^+/G$ is a topological torus again, and so $G$ has rank two. This concludes the proof.
\end{proof}
The following remarkable result by Mess will play a crucial role in the proof of Theorem \ref{th:spacelike}:

\begin{theorem}[Mess \cite{mess}, \cite{gold-mar}] \label{th:mess}
  If $G \subset \mbox{Iso}(\l^3)$ acts freely and properly on $\l^3,$ then $G$ cannot be isomorphic to the fundamental group of a closed  surface of negative Euler characteristic.
\end{theorem}

The following theorem is a consequence of Lemmas \ref{lem:time+light}, \ref{lem:traslaciones} and Theorem \ref{th:mess}.

\begin{theorem}\label{th:spacelike}
Suppose there exists an entire maximal immersion $X:\mb \to \N$ of finite type, where $\N=\l^3/G,$ $G\neq\{\mbox{Id}\}.$ 

Then $\go$ is a group of spacelike translations of rank one or two, and if $\N_+^\uparrow \neq \N$ (i.e., $G \neq \go$) only one of the following possibilities hold:
\begin{enumerate}[(a)]
\item If $G_+^\uparrow=G^\uparrow \neq G_+=G,$ then either $G=\langle R_0\rangle$ or
$G=\langle R_0,T_0\rangle,$ where  in a suitable $(2,1)$-coordinate system  $R_0((x_1,x_2,x_3)):=(x_1+\nu,-x_2,-x_3),$ $\nu \neq 0$ and $T_0(x)=x+(0,\lambda,0),$ $\lambda \neq  0.$ 
\item If $G_+^\uparrow=G_+ \neq G^\uparrow=G,$ then either $G=\langle R_1\rangle$ or  $G=\langle R_1,T_1\rangle,$ where  in a suitable $(2,1)$-coordinate system  $R_1((x_1,x_2,x_3))=(x_1+\delta,-x_2,x_3)$ and $T_1(x)=x+(0,\lambda,0),$ $\delta \neq 0,$ $\lambda \neq 0.$ 

\item If $\go=G_+=G^\uparrow \neq G,$ then either $G=\langle R_2\rangle$ or  
$G=\langle R_2,T_2\rangle,$ where in a suitable $(2,1)$-coordinate system  $R_2((x_1,x_2,x_3))=(x_1,x_2+\delta,-x_3)$ and $T_2(x)=(\lambda,\mu,0),$  $\delta \neq 0,$ $\lambda \neq 0.$ 
\item If  $G_+^\uparrow \neq G_+ \neq G$ and $G_+^\uparrow \neq G^\uparrow \neq G,$ then $G=\langle R_0,R_2\rangle,$ where in a suitable $(2,1)$-coordinate system $R_0$ and $R_2$ are as above. 
\end{enumerate}
Moreover, $X$ is complete and proper,  and $\mb$ is compact if and only if $G_+^\uparrow$ has rank two. If $G=\go$ and $X_0:\mb_0 \to \l^3$ is a conformal parameterization of $X,$  then $\mb_0$ is of finite conformal type and the Weierstrass data of $X_0$ extend  meromorphically to the ends.
\end {theorem}
\begin{proof} As usual,  call  $\tilde{X}:\tilde{\mb} \to \l^3$ a lifting of $X.$
Since $\frac{\Pi_1(\mb)}{\mbox{Ker}(X_*)} \stackrel{X_*}{\cong} G$ (see Remark \ref{re:grupo}),   the surface $\tilde{\sb}:=\tilde{X}(\tilde{\mb})$ is invariant under $G,$ $G$ acts freely and properly on $\tilde{\mb}$ as group of isometries  and $\mb_+^\uparrow:=\tilde{\mb}/\go$ is a finitely sheeted covering of $\mb\equiv \tilde{\mb}/G.$ Moreover,  $X$ induces an entire maximal immersion of finite type $Y:\mb_+^\uparrow \to \N_+^\uparrow,$ where $\N_+^\uparrow$ is the Lorentzian manifold $\l^3/{\go}.$ 

We have to show that $\go$ is a translational group of rank at most $2.$

Reasoning by contradiction suppose that $\go$ is not translational. Lemmae \ref{lem:time+light} and \ref{lem:traslaciones} give that $\go$ consists of orthochronous hyperbolic screw motions, besides the identity. 

Let us see that $\mb_+^\uparrow$ is compact.  Suppose the contrary, and take an annular end $A_0$ of $\mb_+^\uparrow$  containing no singular point. By Theorem \ref{th:parabo}, $A_0$ is is conformally equivalent to a once punctured disc, and any curve $\gamma_0$  generating $\Pi_1(A_0)$ satisfies $\gamma_0 \notin \mbox{Ker} (Y_*).$

Call $R_0\in \go-\{\mbox{Id}\}$ the transformation such that $\tilde{A}_0/\langle R_0\rangle\equiv A_0,$ where $\tilde{A}_0$ is a connected component of $\pg^{-1}(A_0) \subset \tilde{\mb}$ and  $\pg:\tilde{\mb} \to \mb_+^\uparrow$ is the covering projection. Let $N_0:\tilde{A}_0 \to\h^2_-$ denote the Gauss map of $\tilde{X}|_{\tilde{A}_0}.$
Since $N_0\circ R_0=\vec{R}_0\circ N_0,$ we can naturally induce a holomorphic map
$B:A_0\to \h^2_-/\langle \vec{R}_0\rangle$ such that $$B_*:\Pi_1(A_0)\equiv \langle R_0\rangle \longrightarrow\Pi_1(\h^2_-/\langle \vec{R}_0\rangle)\equiv \langle \vec{R}_0\rangle$$ is a group isomorphism.
On the other hand, since $\vec{R}_0$ is an hyperbolic transformation, the quotient $\h^2_-/\langle \vec{R}_0\rangle$ endowed with the hyperbolic metric is conformally equivalent to an annulus $A_1:=\{z\in \C \;:\; 1<|z|<r\},$ $r>1,$ and up to natural identifications, $B:A_0\to A_1$ is  conformal. By the Riemann removable singularity theorem, $B$ extend analytically to the puncture of $A_0,$ and so, $B(\gamma_0)$ is homotopically trivial. This is absurd and proves that $\mb_+^\uparrow,$ and so $\mb,$  must be compact. 

Consider the entire pseudo-spacelike graph $\tilde{\sb}^+$ given in Definition \ref{def:s+}. Since $G_+^\uparrow$ leaves $\tilde{\sb}^+$ invariant,  $\sb_+^\uparrow:=\tilde{\sb}^+/G_+^\uparrow \subset \N_+^\uparrow$ is an embedded incompressible topological surface in $\N_+^\uparrow.$ 
Moreover, as $\sb_+^\uparrow$ is a closed subset of $Y(\mb_+^\uparrow),$  it is compact. Therefore, $G_+^\uparrow$ is cocompact, that is to say, it is isomorphic to the fundamental group of a compact surface (namely, $\Pi_1(\sb_+^\uparrow)$), which contradicts Mess Theorem and proves that $\go$ is translational. Furthermore, by Lemma \ref{lem:traslaciones},  $\go$ is generated by either one or two translations.

For the classification part of the theorem, it is important to keep in mind that an isometry of $G$ reverses the orientation of the spacelike graph $\tilde{\sb}^+$ if and only if it is either positive and non orthochronus or negative and orthochronous. In addition, $\tilde{S}^+/G$ is either a cylinder, a Möbius strip, a torus or a Klein Bottle, and so $G\cong \Pi_1(\tilde{S}^+/G)$ is isomorphic to $\z,$ $\z \times \z$ or $F(a,b)/\langle aabb \rangle,$ where $F(a,b)$ is the non Abelian free group of rank two.

Suppose $\N$ is orientable, i.e., $G_+=G$ (and so, $\go=G^\uparrow \neq G$).

If $\go$ is generated by a spacelike translation, then $G$ is generated by a non orthochronous (hyperbolic, by Lemma \ref{lem:time+light}
) screw motion $R_0.$ Since $R_0^2\in \go,$ the  angle of $R_0$ must be zero, and  in a suitable $(2,1)$-coordinate system $R_0$ is given as in $(a).$ 
When $\go$ has rank two, $G=\langle T'_0,R_0\rangle,$ where $T'_0$ is a spacelike translation, $R_0$ is a non orthochronous hyperbolic screw motion and $\go=\langle R_0^2,T'_0\rangle$ (note that $\go$ is an index two subgroup of $G \cong \Pi_1(\tilde{S}^+/G)$ and $\tilde{S}^+/G$ is a Klein bottle). Reasoning as in the previous case, in a suitable $(2,1)-$coordinate system $R_0((x_1,x_2,x_3))=(x_1,-x_2,-x_3)+(\nu,0,0),$ $\nu \neq 0,$ .  Label $v=(\gamma,\lambda,\mu)$ as the translation vector of $T'_0.$ Since $T_0'':=R_0^{-1}\circ T'_0\circ R_0$ is a translation of vector $(\gamma,-\lambda,-\mu)$ and $\langle T_0'' \circ T_0',R_0^2\rangle$ acts freely and properly on the $x_1$-axis, we get $\gamma=n \delta,$ $n\in\z.$ 
As $T'_0\circ R_0^{-n}$ has no fixed points we infer that $n$ is even. Therefore $T_0:=T'_0\circ R_0^{-n}$ is a  spacelike translation of vector orthogonal to $(\nu,0,0)$ and $G=\langle T_0,R_0\rangle.$ Without loss of generality,  $T_0(x)=(0,\lambda,0),$ $\lambda \neq 0,$  completing $(a).$

In the following we suppose that $G_+ \neq G$.

If $G$ is cyclic, Lemma \ref{lem:time+light} gives $G=\langle R\rangle,$ where, in a suitable $(2,1)$-coordinate system, either $R=R_1$ (orthochronous case, first case of $(b)$) or $R=R_2$ (non orthochronous case, first case of $(c)$). 

Assume that $\go=\gp$ and $G$ has rank two. In this case $G$ is generated by a spacelike translation $T'$ and a negative isometry $R$ without fixed points, where $\go=\langle R^2,T'\rangle.$ As above and in a suitable $(2,1)$-coordinate system, either $R=R_1$ (orthochronous case) or $R=R_2$ (non orthochronous case). Write $T'(x)=x+(v_1,v_2,v_3).$ 
In case $R=R_1,$ $T''=R_1 \circ T' \circ R_1^{-1}$ is a translation of vector $(v_1,-v_2,v_3)$ different from $T'$ ($G$ is the fundamental group of a Klein Bottle) and so $v_2\neq 0.$  Since $T'' \in \langle R_1^2,T'\rangle$ and this group only contains translations of spacelike type, it is not hard to see that $v_3=0,$ and so $T' \circ T''(x)=x+(2v_1,0,0).$ As $\langle T' \circ T'',R_1^2\rangle$ acts freely and properly on $\l^3,$ we get $v_1=k \delta,$ $k \in \z.$ Moreover, since  $T''\circ R_1^{-k}$ has no fixed points, $k$ is even. Thus, $T_1:=T'' \circ R_1^{-k}$ is the translation of vector $(0,\lambda,0),$ where $\lambda=v_2,$ and $G=\langle R_1,T_1\rangle,$ which corresponds to the second case in $(b).$
In case $R=R_2,$ $G$ is commutative and so $T' \circ R_2=R_2 \circ T'.$ This implies that $T_2=T'$ is a horizontal translation, leading to the second case of $(c).$ 

It remains to study the case when $G$ has rank two and $\go \neq \gp.$ In this case $G$ is the fundamental group of a Klein Bottle and  $G/\go \cong \z_2 \times \z_2.$ Therefore,  $G=\langle R_0,R\rangle,$ where $R_0$ is a positive non orthochronous isometry, $R$ is a negative orthochronous isometry and $R_0^2 \circ R^2=\mbox{Id}.$ Since $R_0^2\in \go,$ in a suitable $(2,1)-$coordinate system we have $R_0(x_1,x_2,x_3)=(x_1+\nu,-x_2,-x_3),$ $\nu \neq 0.$ 
As $\vec{G}:=\{\vec{R} \,:\, R \in G\} \cong \z_2 \times \z_2,$ then $\vec{R}^2=\mbox{Id}$ and $\vec{R}$ and $\vec{R}_0$ commute. Thus we can deduce that either  $\vec{R}((x_1,x_2,x_3))=(x_1,-x_2,x_3)$ or  $\vec{R}((x_1,x_2,x_3))=(-x_1,x_2,x_3).$ The last case is impossible because $R \circ R_0$ has no fixed points. Hence $R((x_1,x_2,x_3))=(x_1,-x_2,x_3)+(v_1,v_2,v_3),$ and using that $R_0^2 \circ R^2=\mbox{Id}$ we get $v_1=-\nu$ and $v_3=0.$ Defining $R_2=R_0 \circ R,$ we get  $G=\langle R_0,R_2\rangle,$ which corresponds to $(d).$\\

To finish the theorem, use Theorem \ref{th:parabo} to infer that $Y$ (and so $X$ and $\tilde{X}$) is proper, and that  any annular end of $\mb_+^\uparrow$ (and so, of $\mb$) is conformally equivalent to $\d-\{0\}.$ Since $\tilde{X}$  is a finitely sheeted multigraph and $\go$ consists of translations, $Y(\mb_+^\uparrow)$ (and so $\mb_+^\uparrow$ and $\mb$) is compact if and only if $\tilde {\sb}^+/\go$ and  $\tilde {\sb}^-/\go$ are compact. However, $\tilde{\sb}^+$ and $\tilde{\sb}^-$ are entire graphs, and the last only occurs when $\go$ has rank $2.$

The completeness of $X$ is clear in the compact case. Hence, suppose that $\go$ has rank $1,$ and let us see that $Y$ is complete. 
Let $(g,\phi_3)$ denote the Weierstrass data of $Y_0,$ where $Y_0$ is a conformal parameterization of $Y.$ By the Riemann's removable singularity theorem,  $g$ extends holomorphically to the puncture of any annular end $A_0$ containing no singular point (recall that $|g|<1$), and so $|g||_{A_0}\leq 1-\epsilon,$ $\epsilon>0.$ In the sequel, and up to a Lorentzian isometry, we assume that $\go$ consists of horizontal translations, and so isometrically $\l^3/\go \equiv \Sigma \times \r,$ where $\Sigma$ is a spacelike flat cylinder and $\r$ represents its timelike orthogonal direction.  
If $ds_0^2$ denotes the induced Euclidean metric in $\Sigma$ and $\pi:\l^3/\go \to \Sigma$ the natural projection, then 
$(\pi\circ Y)^*(ds_0^2)\leq ds_1^2 \leq|\frac{\phi_3}{g}|^2,$ where $ds_1^2 = |\phi_1|^2 + |\phi_2|^2 + |\phi_3|^2 = \frac{1}{4} |\phi_3|^2(|g| + 1/|g|)^2.$We infer that the flat metric $|\frac{\phi_3}{g}|^2 $ is complete, and so the same holds for $ds^2|_{A_0} = \frac{1}{4} |\phi_3|^2 (|g| - 1/|g|)^2 \geq    C |\frac{\phi_3}{g}|^2,$ for suitable $C>0.$ Moreover, the completeness of $|\frac{\phi_3}{g}|^2 $ implies that  $\frac{\phi_3}{g}$ (and so $\phi_3$)  extends meromorphically to the ends (see Osserman \cite{osserman}).

\end{proof}

\begin{corollary}\label{co:embedded}
Any entire embedded maximal surface $\sb$ with a finite number of singularities in a complete flat Lorentzian manifold $\N=\l^3/G,$ $G \neq \{\mbox{Id}\},$ is incompressible, proper, complete and of finite type.

Moreover,  up to isometries, only the following  possibilities can hold:
\begin{enumerate}[(a)]
\item $\N=\l^3/\langle  T \rangle,$ where $T$ is a spacelike translation and $\sb$ is a cylinder of finite conformal type.
\item  $\N=\l^3/\langle  T_1,T_2 \rangle,$ where $T_1,$ $T_2$ are independent spacelike translations and $\sb$ is a torus.

\item $\N=\l^3/\langle  R_0 \rangle,$  where $R_0((x_1,x_2,x_3)):=(x_1+\nu,-x_2,-x_3),$ $\nu \neq 0,$ and  $\sb$ is a Möbius strip of finite conformal type.
\item $\N=\l^3/\langle  R_0,T_0 \rangle,$ where $T_0(x)=x+(0,\lambda,0),$ $\lambda \neq  0,$ $R_0$ is as above and $\sb$ is a Klein Bottle.
\item $\N=\l^3/\langle  R_1 \rangle,$ where $R_1((x_1,x_2,x_3))=(x_1+\delta,-x_2,x_3),$ $\delta \neq 0,$ and $\sb$ is a Möbius strip of finite conformal type.
\item  $\N=\l^3/\langle  R_1,T_1 \rangle,$ $R_1$ is as above and $T_1(x)=x+(0,\lambda,0),$ $\lambda \neq 0,$ and $\sb$ is a Klein Bottle.
\item  $\N=\l^3/\langle  R_2 \rangle,$ where $R_2((x_1,x_2,x_3))=(x_1,x_2+\delta,-x_3),$ $\delta \neq 0,$ and $\sb$ is a cylinder of finite conformal type.
\item $\N=\l^3/\langle  R_2,T_2 \rangle,$  where $R_2$ is as above and $T_2(x)=(\lambda,\mu,0),$ $\lambda \neq 0,$ and  $\sb$ is a torus. 
\item $\N=\l^3/\langle  R_0,R_2 \rangle,$ where  $R_0$ and  $R_2$ are as above and $\sb$ is a Klein Bottle.
\end{enumerate}

\end {corollary}
\begin{proof} Let $X:\mb \to \l^3$ be a parameterization of $\sb,$ and write $\tilde{X}:\tilde{\mb} \to \l^3$ its corresponding lifting.

By Theorem \ref{th:spacelike}, $X$ is proper, complete and any annular end of $\mb$ is of finite conformal type.
Since  $\tilde{X}$ is a complete embedding , Lemma \ref{lem:graph} give that $\tilde{\sb}:=\tilde{X}(\tilde{\mb})$ is an entire spacelike graph with isolated singularities over any spacelike plane. This implies that $X$ is imcompressible and  of finite type.

The cases listed above follow Theorem \ref{th:spacelike} again.
\end{proof}

\begin{remark}
The hypothesis in Theorem \ref{th:spacelike} of being maximal is crucial. 
As a matter of fact,  any globally hyperbolic complete Lorentzian manifold admits a smooth Cauchy hypersurface, (see \cite{miguel} for details). 
\end{remark}

\subsection{Maximal surfaces of finite type in $\l^3/\langle  T \rangle$}

Throughout this subsection,  $X:\mb \to \l^3/\langle  T \rangle$ will be an entire maximal surface of finite type, where $T$ is a non trivial spacelike translation. We write $\tilde{X}:\tilde{\mb} \to \l^3$ its corresponding singly periodic lifting, and denote by $(g,\phi_3)$ the Weierstrass data of $X_0:\mb_0 \to \l^3/\langle  T \rangle,$ where $X_0$ is a conformal parameterization of $X.$

We start with the following lemma:
\begin{lemma} \label{lem:final}
The 1-form $\omega=\frac{\phi_3}{g}$ has simple poles at the ends of $\mb_0.$  In particular,  $X$  is asymptotic at any end of $\mb_0$ to a totally geodesic spacelike half cylinder in $\l^3/\langle  T \rangle$ (possibly with multiplicity).
\end{lemma}
\begin{proof} By Theorem \ref{th:spacelike}, $\mb_0$ is of finite conformal type and the Weierstrass data of $X_0$ extend meromorphically to the ends of $\mb_0.$ 
Let $(D^* \equiv \d-\{0\},z)$ be a once punctured conformal disc in ${\mb}_0$ centered at an end of $\mb_0,$ and up to a Lorentzian isometry, suppose that the tangent plane at this end is horizontal, that is to say, $g(0)=0.$ Write $g(z)=z^p,$ $p\geq 1,$  and put 
$$\omega=\left(\sum_{j=-q}^{+\infty} c_{j}z^j\right)dz, \quad c_{-q}\neq 0,$$ on $D.$ Reasoning by contradiction,  suppose $q\geq 2.$ 

Consider the universal covering $\tilde{D}^* \equiv \{u \in\c\;:\;\mbox{Re}(u)<0\} \subset \tilde{\mb}_0$ of $D^*,$ and write $\pi:\tilde{D}^* \longrightarrow D^*,$ $\pi(u)= e^{u},$ 
the covering projection.

Denote by $\pi_0:\l^3 \to \{x_3=0\} \equiv \c$ the orthogonal projection. Without loss of generality, we can suppose that  $\tilde{X}:\tilde{D}^* \to \l^3$ is an embedding and $(\pi_0 \circ \tilde{X})(\tilde{D}^*)$  is a graph over $\{x_3=0\}$ (see the proof of Theorem \ref{th:parabo}).

Using complex notation and following  equation (\ref{eq:wei}), it is not hard to check that
$$\pi_0 \circ \tilde{X} (u)=\frac{-i\, \overline{c}_{-q}}{2 (1-q)}\, e^{(1-q) \overline{u}}\, h(u) \quad u \in \tilde{D}^*,$$
where $h:\tilde{D}^* \to \c$ is  bounded and differentiable in $\tilde{D}^*$ and $\lim_{u \to \infty} h(u)=1$ on the strips $\{u \in \tilde{D}^* \;:\;|\mbox{Im}
(u)|<C\},$ $C \in \r^+.$

Let ${\cal A}:\tilde{D}^* \to \r$ be a smooth branch of the argument ${\cal A}(u):=\arg{\left(\pi_0(\tilde{X}(u))\right)},$ and observe that $A_\theta:=\lim_{r \to -\infty} A(r +i \theta)=\frac{3\pi}{2}-\arg{(c_{-q})}+(q-1) \theta.$

In particular, the graph $\tilde{S}(k):=\tilde{X}(\{u \in \tilde{D}^* \;:\; \mbox{Im}(u) \in [\frac{2 \pi k}{q-1},\frac{2 \pi (k+1)}{q-1}]\}$ is projected by $\pi_0$ over a planar region in $\{x_3=0\}$ containing, up to a compact subset,  the complement of a sector of arbitrarily small angle bisected by the half line $\{(x_1,x_2,0) \;:\; \arg{(x_1,x_2)}=\frac{3\pi}{2}-\arg{(c_{-q})}\},$ and this for any $k \in \z.$ As $\tilde{\sb}(k) \subset \tilde{X}(\tilde {D}^*)$ for any $k \in \z,$ we contradict that $\tilde{X}(\tilde {D}^*)$ is a graph and prove the first part of the lemma.

For the second one, use the same notation as above for $q=1$ and obtain
$$\tilde{X} (u)=\big(\frac{-i\,\overline{c}_{-1}}{2}\, \overline{u},0\big)+ H(u),$$
where  $H:\tilde{D}^* \to \r^3\equiv\c\times\r$ is a bounded differentiable function in $\tilde{D}^*$ and 
$\lim_{u \to \infty} H(u)=(1,\mu_0)$  on the strips 
$\{u \in \tilde{D}^* \;:\; |\mbox{Im}(u)|<C\},$ $C \in \r^+.$ 
Therefore, the surface $\tilde{X}(\tilde{D}^*)$ is asymptotic, as $\mbox{Re(u)}\to -\infty,$ to a horizontal plane and the end $X(\D^*)$ is asymptotic to a totally geodesic horizontal half cylinder invariant under $T^n,$ for suitable $n \in \N.$
\end{proof}

Suppose $\mb_0$ has $r$ ends and let $\overline{\mb}_0=\mb_0 \cup \{P_1,\ldots,P_r\}$ be the conformal compactification of $\mb_0.$ Let  $A_j \subset \mb_0$ be an annular end around $P_j$ containing no singular points, and call $\tilde{A}_j \subset \tilde{\mb}$ a  lifting of $A_j.$ By definition, the integer $w_j \in \n$ such that $\tilde{A}_j/\langle  T^{w_j} \rangle \equiv A_j$ is called the {\em multiplicity} of $P_j.$ When $w_j=1,$ we say that $P_j$ is a Sherk's type end. If $\pg:\tilde{\mb} \to \mb$ is the natural covering, it is clear that $\pg^{-1}(A_j)$ consists of $w_j$ connected components corresponding to the different liftings of $A_j.$ Since the images under $\tilde{X}$ of two different liftings of $A_j$ differ by a translation $T^m,$ $0< m<w_j,$ we infer that  $\tilde{X}$ (and so $X$) is not an embedding provided that $w_j>1,$ $j\in\{1,\ldots,r\}.$   

Let $\Sigma$ be a timelike flat cylinder in $\l^3/\langle  T \rangle,$ and label $H^+$ and $H^-$ as the two connected components of  $\l^3/\langle  T \rangle-\Sigma.$ Without loss of generality, we  suppose that $X(A_j)\cap \Sigma=\emptyset,$ for any $j.$ Therefore we can assign to any end $P_j$ the signature $\varepsilon_j=+1$ or $\varepsilon_j=-1$ depending on $X(A_j) \subset H^+$ or $X(A_j) \subset H^-,$ respectively.

On the other hand, Stokes formula gives 
$\sum_{j=1}^r \mbox{Residue}(\Phi,P_j)=0,$ which leads to $$\sum_{j=1}^r \varepsilon_j w_j=0.$$

The geometrical interpretation of the multiplicity of an end is given in the following corollary:

\begin{corollary} \label{co:jorge1} If we denote by $\bar \mb=\mb \cup \{P_1,\ldots,P_r\}$ the natural topological compactification of $\mb$ and consider $\Pi \subset \l^3$ a spacelike plane invariant under $T,$  then the map $\pi \circ X:\mb \to \Pi/\langle  T \rangle$ extends to a  finitely sheeted branched covering $\pi \circ X:\bar \mb \to \big(\Pi/\langle  T \rangle  \cup \{\infty_+,\infty_-\}\big),$ where  $\pi:\l^3/\langle  T \rangle \to \Pi/\langle  T \rangle$ is the orthogonal projection and  $\Pi/\langle  T \rangle  \cup \{\infty_+,\infty_-\}\equiv \bar \c$ is the natural compactification of the cylinder $\Pi/\langle  T \rangle.$ 

Moreover, the multiplicity of the branched covering $\pi \circ X$ at an end $P_j$ coincides with $w_j.$ Therefore, if up to relabeling we suppose $\{P_1,\ldots,P_s\},$ $s <r,$ are the ends with positive signature, then $\pi \circ X$ is a branched covering of $\sum_{j=1}^s w_j$ sheets.
\end{corollary}

\begin{corollary} \label{co:convexhull1}
The immersion $X$ is an embedding if and only if $r=2$ and $w_1=w_2=1.$ In this case, the convex hull of the graph $\tilde{X}(\tilde{\mb})$ is either a wedge $W$ (if the two Sherk ends of $X$ are not parallel) or a slab (if $X$ has parallel ends). 
\end{corollary}

\subsection{Maximal surfaces of finite type in $\l^3/\langle  T_1,T_2 \rangle$}
Let $X:\mb \to \l^3/\langle  T_1,T_2 \rangle$ be an entire maximal surface of finiten type, where $T_1,$ $T_2$ are independent translations leaving invariant a spacelike plane $\Pi,$  and write $\tilde{X}:\tilde{\mb} \to \l^3$ its corresponding doubly periodic lifting. Label  $\pi:\l^3/\langle  T_1,T_2 \rangle \to \Pi/\langle  T_1,T_2 \rangle$ as the orthogonal projection, and call $(g,\phi_3)$ the Weierstrass data of $X_0:\mb_0 \to \l^3/\langle  T_1,T_2 \rangle,$ where $X_0$ is a conformal parameterization of $X.$ The following corollary is obvious.

\begin{corollary} \label{co:jorge2} The map $\pi \circ X:\mb \to \Pi/\langle  T_1,T_2 \rangle$ is a  finitely sheeted branched covering. 

Moreover, $X$ is an embedding if and only if $\pi \circ X$ is a homeomorphism.
\end{corollary}

\subsection{The case $G=\{\mbox{Id}\}:$ maximal surfaces of finite type in $\l^3$}

A first version of the following theorem was proved in \cite{f-l-s} for embedded maximal surfaces. It shows that the notions of completeness, properness and entireness are equivalent for maximal surfaces in $\l^3$ of finite type (i.e., with a finite number of singularities). We simply sketch the proof.
\begin{theorem} \label{th:compro}
Let $X:\mb\to\l^3$ be  a maximal immersion with  a finite number of isolated singularities. Then the following statements are equivalent:
\begin{enumerate}[(a)]
\item $X$ is complete.
\item $X$ is of finite type
\item $X$ is proper.
\item $X$ is entire.
\end{enumerate} 
In any case, if $X_0:\mb_0 \to \l^3$ is a conformal parameterization of $X,$ $\mb_0$ is biholomorphic to a $\overline{\mb}_0-\{P_1,\ldots,P_r\},$ where $\overline{\mb}_0$ is a compact Riemann surface with analytical boundary and $\{P_1,\ldots,P_r\} \subset \overline{\mb}_0-\partial (\overline{\mb}_0).$ Moreover, the Weierstrass data $(g,\phi_3)$ of $X_0$ extend meromorphically to $\overline{\mb}_0.$ 
\end{theorem}
\begin{proof} $(a)\Longrightarrow (b):$ Huber's Theorem \cite{huber} implies that $\mb_0$ has finite conformal type, that is to say, it has finite topology and its ends are conformal once punctured disc. Osserman's \cite{osserman} classical results  imply that the Weierstrass data $(g,\phi_3)$ extend meromorphically to the ends, and the pole order of $\Phi:=(\phi_1,\phi_2,\phi_3)$ at the ends is at least $2.$ Finally, one can use Jorge-Meeks \cite{jorgemeeks} ideas to control the asymptotic behavior at infinity, proving in particular that $\pi \circ X:\mb \to \Pi$ is a finitely sheeted covering, where  $\pi$ is the orthogonal projection over any spacelike plane.

$(b)\Longrightarrow (c)$ is trivial, and $(c) \Longrightarrow (d)$ follows from Lemma \ref{lem:graph}.

To check $(d)\Longrightarrow (a),$ let $K_0 \subset \mb$ be a compact subset containing the singularities. We can assume without loss of generality that $|g|<1$ on $\mb-K_0,$ because$X|_{\mb-K_0}$ is spacelike.On the other hand, if $ds_0^2$ is the induced Riemannian metric on $\Pi,$ we observe that $(\pi\circ X)^*(ds_0^2)\leq ds_1^2 \leq|\frac{\phi_3}{g}|^2,$ where$ds_1^2 = |\phi_1|^2 + |\phi_2|^2 + |\phi_3|^2 = \frac{1}{4} |\phi_3|^2(|g| + 1/|g|)^2.$Thus, the flat metric $|\frac{\phi_3}{g}|^2 $ is complete, and so itfollows from  classicalresults of Huber \cite{huber} and  Osserman \cite{osserman} that$\mb-K_0$  is conformally a finitely punctured compact Riemann surface with compact boundary, $\phi_3/g$ has  poles at the puncturesand $g$ extends holomorphically to the punctures.Consequently, $|g|<1-\epsilon$ on $\mb-K_0$  for some$\epsilon >0,$  and so $ds^2 = \frac{1}{4} |\phi_3|^2 (|g| - 1/|g|)^2 \geq\frac{\epsilon^2}{4} |\frac{\phi_3}{g}|^2$ is complete.
\end{proof}

\begin{corollary} \label{co:jorge} 
Let $X:\mb \to \l^3$ be a entire maximal immersion of finite type. With the notation of Theorem  \ref{th:compro}, label $\overline{\mb}=\mb \cup \{P_1,\ldots,P_r\}$ as the natural  compactification of $\mb.$ Let  $\Pi$ be a spacelike plane and call  $\pi:\l^3 \to \Pi$ the orthogonal projection. 

Then the map $\pi \circ X:\mb \to \Pi$ extends to a  finitely sheeted branched covering $\pi \circ X:\overline{\mb} \to \Pi\cup \{\infty\}\equiv \bar \c.$  Moreover, the end $P_j$ is a branch point of $\pi \circ X$ with multiplicity
  $$w_j:= \mbox{Ord}_{P_j}(\Phi) -1 = \max \{\mbox{Ord}_{P_j}(\phi_k),\;k=1,2,3\}-1 \geq 1,$$  where $(\phi_1,\phi_2,\phi_3)$ is the Weierstrass representation of $X$ and $\mbox{Ord}_{P_j}(\phi_k)$ is the pole order of $\phi_k$ at $P_j,$ $j=1,\ldots,r,$ $k=1,2,3.$ 
In particular, $\pi \circ X$ is a branched covering with $\sum_{j=1}^r w_j$ sheets.

As a consequence, $X$ is an embedding if and only if $r=1$ and $w_1=1.$
\end{corollary}
\begin{proof}
Let $A$ be any annular end in $\mb$ corresponding to $P_j,$ and put $h:=\pi \circ {X}.$

Let $U \subset \{x_3=0\}$ be an open disc containing $h(\partial(A))$ and define $A'$ as the non compact connected component of $A \cap h^{-1} (\{x_3=0\}-U).$ Obviously, $h(A') \cap (\{x_3=0\}-U)$ is open and closed in $\{x_3=0\}-U$ (we use here that $h:{\mb} \to \{x_3=0\}$ is a finitely sheeted multigraph), and so  $h(A') = \{x_3=0\}-U$. We infer that $h|_{A'}: A' \to \{x_3=0\}-U$ is a proper local homeomorphism, and so a finitely sheeted covering. Since $\phi_j,$ $j=1,2,3,$ have no real periods, it is easy to check that $w_j \geq 1.$ Straightforward arguments give that $h|_{A'}$ has $w_j \geq 1$ sheets (see \cite{jorgemeeks} for a similar discussion in  the case of minimal surfaces), which proves that $P_j$ is a branch point of multiplicity $w_j$ and concludes the proof.  
\end{proof}

\section{Construction of entire maximal surfaces of finite type} \label{sec:singly-doubly}

Throughout this section, and following Theorem \ref{th:spacelike}, we  suppose that $G\subset\mbox{Iso}(\l^3)$ is a discrete group of spacelike translations of rank $r(G)=0,$ $1$ or $2$ (rank $0$ means that $G=\{\mbox{Id}\}$), and $\tilde{X}:\tilde{\mb}\to\l^3$ is an entire maximal immersion invariant under $G$ such that the immersion $X:\mb=\tilde{\mb}/G \to \l^3/G$ has finite type, $k_1$ lightlike singularities, $k_2$ spacelike singularities and $r$ ends. We write $k:=k_1+k_2$ as the total number of singular points of $X.$

By Theorem \ref{th:parabo}, the conformal support $\mb_0$ of $X$ is biholomorphic to $\rb- \big( (\cup_{j=1}^{k_1} D_j) \cup \{P_1,\ldots,P_r\}\big),$ where $\rb$ is a compact Riemann surface of  genus $\xi_0:=\mbox{Gen}(\mb)=\mbox{Gen}({\mb}_0),$ $\{D_j \;:\; j=1 \ldots,k_1\}$ are open conformal discs in $\rb$ with pairwise disjoint clusures and $\{P_1,\ldots,P_r\} \subset \rb- \cup_{j=1}^{k_1} \bar D_j.$ 

We denote by $(g,\phi_3)$ the meromorphic data associated to the conformal parameterization $X_0:\mb_0 \to \l^3/G$ of $X.$ Up to a Lorentzian isometry, we always assume that $|g| \leq 1$ on $\mb_0.$ 

If $r(G)<2$ (non compact case),  $(g,\phi_3)$ extend meromorphically to the ends $\{P_1,\ldots,P_r\}$ of $\mb_0$ and the induced metric on $\mb_0,$ namely $ds^2=\frac{1}{4}(1-|g|^2)^2|\frac{\phi_3}{g}|^2,$ is complete (see Theorems \ref{th:spacelike} and \ref{th:compro}). Corollary \ref{co:jorge} and Lemma \ref{lem:final} imply that the 1-form $\omega:=\frac{\phi_3}{g}$ has poles at the ends of order $\geq 2$ (if $r(G)=0$) or equal to $1$ (if $r(G)=1$). This simply means that $\omega$ is meromorphic at the compact surface $\rb-\cup_{j=1}^k D_j.$ In the doubly periodic case $\mb_0$ is compact and $\omega$ holomorphic. 

Write  $\partial (\mb_0) =\cup_{j=1}^{k_1} \gamma_j,$ $\gamma_j:=\partial (D_j) \equiv \s^1,$ and let $\mb_0^*$ be the mirror of $\mb_0$ with boundary  $\partial (\mb_0^*) =\cup_{j=1}^{k_1} \gamma_j^*,$ $\gamma_j^* \equiv \s^1$ (see Lemma \ref{lem:sing1}). 
Let $\Sg\equiv\mb_0 \cup \mb_0^*$ the double of $\mb_0,$ that is to say, the Riemann surface without boundary obtained by gluing analytically $\mb_0$ and $\mb_0^*$ along $\gamma_j$ and $\gamma_j^*,$ $j=1,\ldots,k_1.$ The natural mirror involution $J:\Sg \to \Sg,$ $J(p):=p^*,$ fixes pointwise $\gamma_j\equiv \gamma_j^* \subset \Sg,$ $j=1,\ldots,k_1.$  
We denote by $\overline{\Sg}$  the conformal compactification of $\Sg,$ and keep calling $J$ the antiholomorphic mirror involution on $\overline{\Sg}$ (recall that in the doubly periodic case  $\overline{\Sg}=\Sg$). Since $g$ extends meromorphically to $\overline{\Sg},$ it has well defined degree, namely $\mbox{Deg}(g).$

Note that $\overline{\Sg}$ is a compact Riemann surface of genus $2\xi_0+k_1-1.$ From Lemma \ref{lem:sing1} the Weierstrass data $(g,\phi_3)$ of $X_0$ extend, with the same name, to meromorphic data on $\overline{\Sg}$ satisfying $g \circ J=1/\bar g$ and $J^*(\phi_3)=-\overline{\phi}_3,$ that is to say, 
$J^*(\phi_j)=-\overline{\phi}_j,$ $j=1,2,3.$

For convenience, we refer to $(\Sg,J,g,\phi_3)$ as the Weierstrass representation of $X_0.$

In the embedded case, Lemma \ref{lem:graph} and Corollary \ref{co:embedded}  show that either $\mb$ is topologically a plane and $\xi_0=0$ ($r(G)=0$), a cylinder and $\xi_0=0$ ($r(G)=1$) or a torus and $\xi_0=1$ ($r(G)=2$). Furthermore, Lemma \ref{lem:sing1} implies that the singularities are of conelike type.\\

The relationship between the topology of $\mb,$ the degree of $g$ and the geometry of $X$ at the ends is explained in the following theorem. 

Let $Q_1, \ldots, Q_{k_2} \in \mb_0$ denote the spacelike singularities of $X,$ and call $n_j$ as the zero order of $\Phi$ at $Q_j,$ $j=1,\ldots,k_2.$
We also denote by $n_j'$ the total vanishing order of $\Phi$ at $\gamma_j,$ that is to say, the number of zeroes counted with multiplicity of $\Phi$ on $\gamma_j$ (which is always even), $j=1,\ldots k_1.$ Likewise, $m_j$ will denote the degree of $g|_{\gamma_j}:\gamma_j \to \{|z|=1\},$ $j=1,\ldots,k_1.$ Finally, call $w_j$ the multiplicity of $X$ at the end $P_j,$  $j=1,\ldots,r,$ and define $$V_s:=\sum_{j=1}^{k_2} n_j, \quad V_l:=\sum_{j=1}^{k_1} n'_j.$$ 
Since $|g|^{-1}(1)=\cup_{j=1}^{k_1} \gamma_j,$ we get that $\mbox{Deg}(g)=\sum_{j=1}^{k_1} m_j.$ 

Write $W_\infty= \sum_{j=1}^r \mbox{Ord} (\Phi,P_j),$ where $\mbox{Ord} (\Phi,P_j)$ is the pole order of $\Phi$ at $P_j.$ Note that $W_\infty $ is equal to $\sum_{j=1}^r (w_j+1),$ $r$ or $0$ provided that $r(G)=0,$ $1$ or $2,$ respectively.

\begin{theorem}[Topological formula]
With th previous notation, the following formula holds: 
$$k_1-\chi (\overline{\mb})=V_s+\frac{V_l}{2}+\mbox{Deg}(g)-W_\infty,$$ 
where $\overline{\mb}:=\mb \cup \{P_1,\ldots,P_r\}$ is the natural compactification of $\mb$ and
$\chi (\overline{\mb}):=2-2 \xi_0$  the  Euler characteristic of $\overline{\mb}.$
\end{theorem}

\begin{proof}
Let $\Pi$ be a spacelike plane invariant under $G,$ and call $\pi:\l^3/G \to \Pi/G$ the Lorentzian orthogonal projection. Let $(\Pi/G)^*$  be either $\Pi \cup \{\infty\}\equiv \s^2$ (if $G=\{\mbox{Id}\}$), $\Pi/G \cup \{-\infty,+\infty\}\equiv \s^2$ (if $G=\langle  T\rangle$) or $\Pi/G \equiv \s^1\times \s^1$ (if $G=\langle  T_1,T_2 \rangle)$).

Consider the finitely sheeted branched topological covering  $h:=\pi \circ X:\overline{\mb} \to (\Pi/G)^*,$ and use Riemann-Hurwitz formula to get: 
$$\chi(\overline{\mb})=\chi\big((\Pi/G)^*\big) \mbox{Deg}(h)- \pmb{B},$$
where $\pmb{B}$ is the total branching number of  $h.$ By Lemmae \ref{lem:sing} and \ref{lem:sing1} and Corollaries \ref{co:jorge} and \ref{co:jorge1} we get $\pmb{B}=B_s+B_l+B_\infty,$ where $B_s=V_s$ and $B_l=\mbox{Deg}(g)+\frac{V_l}{2}-k_1$ 
are the total branching numbers at the spacelike and lightlike singularities resp., 
and $B_\infty$ is the total branching number at the ends, that is, $B_\infty=\sum_{j=1}^r  (w_j-1)$ if $r(G)<2,$ and $B_\infty =0$ when $r(G)=2.$ 

Furthermore, when $r(G)<2,$ Corollaries \ref{co:jorge} and \ref{co:jorge1} give $\chi\big((\Pi/G)^*\big) \mbox{Deg}(h)=B_\infty + W_\infty.$ This equation trivially holds when $r(G)=2,$ and so in general
\begin{equation*}
\chi(\overline{\mb})=W_\infty -V_s-\frac{V_l}{2}-\mbox{Deg}(g)+k_1
\end{equation*}
This completes the proof.
\end{proof}

The following theorem summarizes all the known information about the Weierstrass representation of entire maximal surfaces of finite type, and provides an useful construction method.

\begin{theorem} [Analytical representation of entire maximal surfaces of finite type] \label{th:periodic}
Let $X:\mb \to \l^3/G$ be an entire maximal immersion of finite type with $k\geq 1$ singularities, where $G$ is a rank $\leq 2$ discrete group of spacelike translations acting freely and properly in $\l^3.$ Let $X_0:\mb_0 \to \l^3/G$ be a conformal parameterization of $X,$ and call $(\Sg:=\mb_0 \cup \mb_0^*,J,g,\phi_3)$ its Weierstrass data. Then:
\begin{enumerate}[(i)]
\item $g\circ J=1/\bar{g}$ and $|g|<1$ on $\mb_0-\partial(\mb_0).$ 
\item The vectorial 1-form $\Phi:=(\phi_1,\phi_2,\phi_3)$ given in equation (\ref{eq:wei}) is holomorphic in $\Sg,$ have no zeroes in $\partial (\mb_0)$  and satisfies $J^*(\Phi)=-\overline{\Phi}.$ 
\item The translations in $G$ are given by the vectors
$\{\mbox{Re}\int_\gamma(\phi_1,\phi_2,\phi_3) \;:\;\gamma \in {\cal H}_1(\mb,\z)\}.$
\item  $r(G)<2$ if and only if $\Sg$ is not compact. Furthermore,  if $r(G)=0$ then $\Phi$ has poles of order $\geq 2$ at the ends of $\Sg,$ and if $r(G)=1$ then $\Phi$ has simple poles at the ends of $\Sg.$ In the compact case ($r(G)=2$),  $\Phi$ is holomorphic.
\end{enumerate}

Conversely, let $\Sg$   and $J:{\Sg}\to {\Sg}$ be a  Riemann surface of finite conformal type, $\partial (\Sg)=\emptyset,$ and an  antiholomorphic involution such that the fixed point set of $J$ consists of $k_1$ pairwise disjoint analytical Jordan curves $\gamma_j,$ $j=1,\ldots,k_1,$ spliting $\Sg$ into two connected regions, whose closures in $\Sg$ will be denoted by $\mb_0$ and $\mb_0^*:=J(\mb_0).$    Assume that $(g,\phi_3)$ are Weierstrass data on $\Sg$ satisfying
$(i),$ $(ii),$ and defining $G$ like in $(iii),$  $r(G) \leq 2$ and $(iv)$ holds.

Then the map 
\begin{equation}\label{eq:inmersion}
X_0:\mb_0 \to \l^3/G,\quad  X_0(p):=\mbox{Re}\int_{p_0}^p \Phi,
\end{equation}
where $p_0\in\mb_0,$ is well defined and induces an entire maximal  immersion of finite type of the surface $\mb$ in $\l^3/G.$ Here, $\mb$ denotes the quotient surface of $\mb_0$ identifying 
each boundary circle $\gamma_j$ with a single point $q_j,$ where $q_n \neq q_m$ provided that $n \neq m$ and $F:=\{q_j \;:\; j=1,\ldots,k\} \cap \mb_0=\emptyset.$ 
  Furthermore, $X$ has $k:=k_1+k_2$ singularities, where $k_2$ is the number of zeroes counted with multiplicity of $\Phi$ in $\mb_0-\cup_{j=1}^{k_1} \gamma_j.$

In addition, $X$ is an embedding if and only if $k_2=0$ and  $g:\gamma_j\to\s^1$ is injective, $j=1,\ldots,k$ (i.e., $\mbox{Deg}(g)=k$). In this case,  $\mb_0$ is conformally equivalent to either $\c$ minus $k$ discs (if $r(G)=0$), $\c-\{0\}$ minus $k$ discs (if $r(G)=1$) or a torus minus $k$ discs (if $r(G)=2$).
\end{theorem}

\begin{proof} The first part of the theorem is just a summary of the previously stablished results.

For the construction part, note that $J_*(\gamma_j)=\gamma_j$ and $(ii)$ give that $\mbox{Re}\left(\int_{\gamma_j} \Phi \right)=0,$ for any $j.$ Therefore, the real periods of $\Phi$ are the vectors in  $G=\{\mbox{Re}\int_\gamma \Phi \;:\;\gamma \in {\cal H}_1(\mb,\z)\}.$

Obviously, the immersion $X_0:\mb_0 \to \l^3/G,$ $X:=\mbox{Re}\int(\phi_1,\phi_2,\phi_3),$ is well defined and leads to an immersion of $\mb$ in $\l^3/G$ with $k$ singularities. By Lemma \ref{lem:graph}, $X$ is an entire spacelike immersion provided that the metric $ds^2$ induced by $X$ on $\mb_0$ is complete. This fact is an elementary consequence of the expresion $ds^2=\frac{1}{4}(1-|g|^2)^2|\frac{\phi_3}{g}|^2$ and $(iv).$ 

It remains to check that $X$ is of finite type. If $r(G)=0,$ $X$ is an entire maximal immersion with a finite number of singularities in $\l^3,$ and so, a finite multigraph by Corollary \ref{co:jorge}. When $r(G)=1,$ the ends of $X$ lift to Scherk's type ends of $\tilde{X},$ as explained in Lemma \ref{lem:final}, and since $X$ has a finite number of ends,  $\tilde{X}$ is a finite multigraph too. The compact case $r(G)=2$  is obvious. 

Finally, note that  $X$ is an embedding if and only if $\tilde{X}$ is, which prevents the existence of spacelike singularities. By Lemmas \ref{lem:graph} and \ref{lem:sing1}, $\tilde{X}$ is an embedding if and only if all its singularities are of conelike type, which concludes the proof.
\end{proof}

\subsection{Examples of entire embedded maximal surfaces of finite type}
To finish, we present some families of entire embedded maximal surfaces of finite type in $\l^3/G,$ $G \neq \{\mbox{Id}\},$  sweeping all the cases in Theorem \ref{th:spacelike}. A large family of embedded examples of fnite type in $\l^3$ can be found in \cite{f-l-s}, and a more general family of non embedded examples can be found in \cite{isabel}.\\

{\bf Examples with one singular point and non parallel ends in $\l^3/\langle  T \rangle.$}

Consider the following data:

$${\Sg}=\overline{\c}-\{b,-b,1/b,-1/b\},\quad J:\Sg\to\Sg,\;J(z)=1/\bar{z}\quad
 \mb_0=\overline{\d(0,1)}-\{b,-b\},\quad b\in]0,1[$$
$$g(z)=z,\quad \phi_3=\frac{z dz}{(z^2-b^2)(b^2 z^2-1)}$$

Then the immersion given by Equation ($\ref{eq:inmersion}$)
is an embedded entire maximal surface of finite type in $\l^3/G,$ where  $G$ is the group generated by the translation $T$ of vector $v=\mbox{Re}(2\pi i Res(\Phi,b))=(\frac{\pi}{2b(b^2+1)},0,0).$ 

Moreover, the antiholomorphic transformation of ${\mb}_0$ given by $A(z):=-\bar{z}$ satisfies 
$A^*(\phi_1,\phi_2,\phi_3)=(\bar{\phi}_1,-\bar{\phi}_2,\bar{\phi}_3).$ 
Thus, the lifted surface in $\l^3,$ $\tilde{X}:\tilde{\mb} \to \l^3,$ is invariant by the isometry 
$$R(x_1,x_2,x_3)=(x_1,-x_2,x_3)+\mbox{Re}\int_0^{A(0)}\Phi=(x_1,-x_2,x_3)$$
and so the quotient of $\tilde{X}$ by the group generated by the isometry $R_1:=T\circ R$ is an embedded entire maximal surface with a singularity corresponding to the first case described in Theorem \ref{th:spacelike} $(b).$ A surface in this family has been illustrated in Figure \ref{fig:scherk}.\\

{\bf Examples with two singular points and parallel ends in $\l^3/\langle  T \rangle.$}

Consider the data:
$${\Sg}=\{(z,w)\in\overline{\c}^2\,:\,w^2=\frac{(z-a)(z-b)}{(az-1)(bz-1)},\; z \neq 0,\infty\}, \quad
J:{\Sg}\to{\Sg},\; J(z,w)=(1/\bar{z},1/\bar{w})$$
$$\mb_0=\{(z,w)\in {\Sg}\,:\,|z|\leq 1\} \quad b<a<1,\,a>0,\,b\neq 0$$
$$g(z,w)=z\quad \phi_3=\frac{dz}{w(az-1)(bz-1)}$$

Then the immersion given by Equation ($\ref{eq:inmersion}$) is an embedded entire maximal surface with two singularities in $\l^3/G,$ where $G$ is generated by the translation $T$ of vector $v=\mbox{Re}(2\pi i Res(\Phi,0))=\mbox{Re}((\frac{-\pi}{\sqrt{ab}},\frac{-i\pi}{\sqrt{ab}},0)).$ Thus, the translation vector is in the direction of the $x_1-$axis in case $b>0$ or in the direction of the $x_2-$axis if $b<0.$

Moreover, if we choose $b=-a,$ the transformations of $\mb_0$ given by 
$A_0(z,w)=(-\bar{z},\bar{w}),$
$A_1(z,w)=(\bar{z},\bar{w})$ and
$A_2(z,w)=(-z,w),$
lift to the following isometries in $\l^3:$
$$R_0(x_1,x_2,x_3)=(-x_1,x_2,-x_3)+\mbox{Re}\int_a^{A_0(a)}\Phi=(-x_1,x_2,-x_3)+\frac{v}{2},$$
$$R_1(x_1,x_2,x_3)=(-x_1,x_2,x_3)\quad {\rm and} \quad R_2(x_1,x_2,x_3)=(x_1,x_2,-x_3)+\mbox{Re}\int_a^{A_2(a)}\Phi=(x_1,x_2,-x_3)+\frac{v}{2}$$

The quotient of the lifted surface $\tilde{X}:\tilde{\mb}\to \l^3$ in $\l^3$ by the groups $\langle R_0\rangle,$ $\langle T\circ R_1\rangle$ and $\langle R_2\rangle$ give examples of embedded entire maximal surfaces  corresponding to the first cases described in Theorem \ref{th:spacelike}  $(a),$ $(b)$ and $(c),$ respectively.

A surface in this family has been illustrated in Figure \ref{fig:periodicas} (right).\\

{\bf Examples with two singular points in $\l^3/\langle  T_1,T_2 \rangle.$}

Consider the data:

$${\Sg}=\{(z,w)\in\overline{\c}^2\,:\,w^2=\frac{(z^2-a_1^2)(z^2-a_2^2)}{(a_1^2z^2-1)(a_2^2z^2-1)}\},\quad 
a_1,a_2\in\r^*,\;a_1\neq a_2$$
$$J:{\Sg}\to{\Sg},\; 
J(z,w)=(1/\bar{z},1/\bar{w}),\quad \mb_0=\{(z,w)\in{\Sg}\,:\,|z|\leq 1\}$$
$$g(z,w)=z, \quad \phi_3=\frac{z dz}{w(a_1^2z^2-1)(a_2^2z^2-1)}$$

Then, Equation (\ref{eq:inmersion}) leads to an embedded entire maximal surface with two singularities in $\l^3/G,$ where $G$ is generated by two horizontal translations $T_1$ and $T_2$ of vectors $v_1$ and $v_2.$ Namely, $v_i=2 \mbox{Re}\int_{\gamma_i}\Phi,$ $i=1,2,$ where $\gamma_1$ (resp. $\gamma_2$) is a simple arc in $\mb_0$ joining $a_1$ and $-a_1,$ (resp. $a_1$ and $a_2$). Elementary computations show that in fact $v_1=(0,\lambda,0),$ and $v_2=(\mu,0,0),$ for suitable $\lambda,\,\mu\in\r^*.$ The surface on the left in Figure \ref{fig:periodicas} is an example of this family of surfaces.

As in the previous cases, the transformations of ${\mb}_0,$ $A_0(z,w):=(\bar{z},-\bar{w}),$ $A_1(z,w):=(\bar{z},\bar{w})$ and $A_2(z,w):=(-z,-w)$ induce isometries of $\l^3,$ $R_0,$ $R_1$ and $R_2$ resp., leaving the lifted immersion $\tilde{X}$ invariant. To be more precise,
$$R_0(x_1,x_2,x_3)=(x_1,-x_2,-x_3), \quad R_1(x_1,x_2,x_3)=(-x_1,x_2,x_3)\quad {\rm and} \quad R_2(x_1,x_2,x_3)=(x_1,x_2-\lambda/2,-x_3)$$
Thus, the corresponding quotients of $\tilde{X}$ by the groups $\langle T_2 \circ R_0, T_1\rangle,$ $\langle T_1 \circ R_1,T_2\rangle,$ $\langle R_2,T_2\rangle$ and $\langle T_2\circ R_0,R_2\rangle$ provide examples of embedded entire maximal surfaces of finite type corresponding to the second cases described in Theorem \ref{th:spacelike} $(a),$ $(b),$ $(c),$ and Theorem \ref{th:spacelike} $(d),$ respectively.

A surface  in this family has been illustrated in Figure \ref{fig:periodicas} (left).


{\bf ISABEL FERNANDEZ, FRANCISCO J. LOPEZ,} \newline
Departamento de Geometr\'{\i}a y Topolog\'{\i}a \newline
Facultad de Ciencias, Universidad de Granada \newline
18071 - GRANADA (SPAIN) \newline
e-mail:(first author) isafer@ugr.es, (second author) fjlopez@ugr.es

\end{document}